\documentclass[pdflatex,sn-mathphys-num]{sn-jnl}% Math and Physical Sciences Numbered Reference Style 
\usepackage{graphicx}%
\usepackage{multirow}%
\usepackage{amsmath}
\usepackage{amssymb}
\usepackage{amsfonts}%
\usepackage{amsthm}%
\usepackage{mathrsfs}%
\usepackage[title]{appendix}%
\usepackage{xcolor}%
\usepackage{textcomp}%
\usepackage{manyfoot}%
\usepackage{booktabs}%
\usepackage{algorithm}%
\usepackage{algorithmicx}%
\usepackage{algpseudocode}%
\usepackage{listings}%
\usepackage{enumitem}
\usepackage{threeparttable}
\usepackage{multirow}
\usepackage{makecell}
\usepackage{adjustbox}
\usepackage{mathrsfs}
\usepackage{amsmath}
\usepackage{natbib}% Citation support using natbib.sty
\usepackage{tikz}
\usepackage{lmodern}

\usepackage{endnotes}
\usepackage{bbm}
\usepackage{pstricks}
\usepackage[justification=centering]{caption}
\usepackage{hyperref}
\hypersetup{
    colorlinks=true,
    linkcolor=blue,
    filecolor=magenta,
    urlcolor=cyan,
}
\usepackage{comment}
\usepackage{mathtools}

\def\cdf(#1)(#2)(#3){0.5*(1+(erf((#1-#2)/(#3*sqrt(2)))))}

\usepackage{pgfplots}
\pgfplotsset{compat=1.17}
\usepgfplotslibrary{fillbetween}

\def\Dp{\mathfrak{p}}

\def\BY{\mathbb{Y}}
\def\XP{\mathcal{X}}
\def\WP{\mathcal{W}}
\def\WprimeP{\mathcal{W'}}
\def\STpolytope{\mathscr{P}}
\def\BE{\mathbb{E}}
\def\Re{\mathbb{R}}
\def\BT{\mathbb{T}}
\def\BN{\mathbb{N}}
\def\BP{\mathbb{P}}
\def\BI{\mathbbm{1}}
\def\Var{\text{VaR}}
\def\CVar{\text{CVaR}}
\def\ie{\text{i.e.,}}
\def\RI{\text{R}_{\text{I}}}
\def\XTP{\mathcal{X}_T}

\def\IQF{\text{I}_{\BT(x)}}

\newtheorem{theorem}{Theorem}%  meant for continuous numbers
\newtheorem{proposition}[theorem]{Proposition}% 
\newtheorem{remark}{Remark}%

\newtheorem{definition}{Definition}%
\newtheorem{lemma}{Lemma}
\newtheorem{corollary}{Corollary}

\begin{document}

\title[Risk Minimization]{Minimizing risk measures with applications in network traffic engineering}

%%=============================================================%%
%% GivenName	-> \fnm{Joergen W.}
%% Particle	-> \spfx{van der} -> surname prefix
%% FamilyName	-> \sur{Ploeg}
%% Suffix	-> \sfx{IV}
%% \author*[1,2]{\fnm{Joergen W.} \spfx{van der} \sur{Ploeg} 
%%  \sfx{IV}}\email{iauthor@gmail.com}
%%=============================================================%%

\author*[1]{\fnm{Ashish} \sur{Chandra}}\email{achand6@ilstu.edu}

\author[2]{\fnm{Mohit} \sur{Tawarmalani}}\email{mtawarma@purdue.edu}

\affil*[1]{\orgdiv{College of Business}, \orgname{Illinois State University}, \orgaddress{\city{Normal}, \postcode{61761}, \state{IL}, \country{USA}}}

\affil[2]{\orgdiv{Mitch Daniels School of Business}, \orgname{Purdue University}, \orgaddress{\city{West Lafayette}, \postcode{47906}, \country{USA}}}

%%==================================%%
%% Sample for unstructured abstract %%
%%==================================%%

\abstract{This paper presents a novel two-stage optimization framework designed to model integrated quantile functions, which leads to the formulation of a bilinear optimization problem (P). A specific instance of this framework offers a new approach to minimizing the Value-at-risk ($\Var$) and the Conditional Value-at-risk ($\CVar$), thus providing a broader perspective on risk assessment and optimization. We investigate various convexification techniques to under- and over-estimate the optimal value of (P), resulting in new and tighter lower- and upper-convex estimators for the $\Var$ minimization problems. Furthermore, we explore the properties and implications of the bilinear optimization problem (P) in connection to the integrated quantile functions. Finally, to illustrate the practical applications of our approach, we present computational comparisons in the context of real-life network traffic engineering problems, demonstrating the effectiveness of our proposed framework.}


%%\pacs[JEL Classification]{D8, H51}

%%\pacs[MSC Classification]{35A01, 65L10, 65L12, 65L20, 65L70}

\maketitle

\section{Introduction}\label{sec1}

Risk measures are essential in decision-making under uncertainty, providing quantifiable insights into potential losses in diverse fields, such as finance, engineering, supply chain management, resource allocation, project planning, etc. These measures help guide effective strategies for managing downside risk, a central concern in many industries. Multiple risk measures exist to address different aspects of these uncertainties; see for instance \cite{dhaene2006risk, emmer2013best, gambrah2014risk}. Among these risk measures, Value-at-Risk ($\Var$) and Conditional Value-at-Risk ($\CVar$) stand out and have been instrumental in providing a quantifiable framework for assessing risk. 
$\Var$, for instance, is widely used for its simplicity and interpretability, providing a clear threshold beyond which a given percentage of worst-case losses may occur. Yet, $\Var$ does not account for the severity of losses beyond this threshold, which can lead to an underestimation of tail risks in high-impact situations. $\CVar$, on the other hand, considers the average loss exceeding $\Var$ and has thus become a preferred measure in fields requiring greater risk sensitivity \cite{rockafellar2000optimization}.
Various quantile-based risk measures \cite{embrechts2018quantile, cont2010robustness} have also been developed to extend traditional $\Var$ and $\CVar$ metrics, addressing their limitations in specific applications. 
The authors in \cite{embrechts2018quantile, cont2010robustness} proposed Range-Value-at-Risk (RVaR) as a general two-parameter family of risk measures that bridges $\Var$ and $\CVar$ by considering an interval of quantiles rather than a single point. \cite{belles2016use} proposed a new family of risk measures GlueVaR which belongs to the class of distortion risk measures and showed that GlueVaR measures can be defined as a linear combination of common quantile-based risk measures ($\Var$, $\CVar$).

In the context of decision-making under uncertainty, risk-averse frameworks have been increasingly applied across various areas, such as portfolio selection, supply chain disruption, machine scheduling, and network traffic engineering, etc. (see, for example, \cite{qazi2023supply, babazadeh2018optimisation, lotfi2020robust, v2001value, gaivoronski2005value, dixit2020project, atakan2017minimizing, sarin2014minimizing, bogle2019teavar, chandra2022techniques}). This widespread application underscores the importance of making decisions that minimize the risk metric under consideration. 
$\Var$ remains one of the most popular risk measures in this context, despite the challenges in its optimization due to its non-smooth and non-convex nature (see for instance \cite{mausser1998beyond, gaivoronski2005value, larsen2002algorithms}). Although optimization methods for $\Var$ have been explored extensively (see, for instance, \cite{gaivoronski2005value, larsen2002algorithms}), efficient algorithms for real-world applications still remain scarce. 
One prominent approach for $\Var$ minimization is the difference-of-convex (DC) function method, which leverages representing $\Var$ as the difference between two convex functions. The authors in \cite{wozabal2010difference, wozabal2012value} use the DC function method for the case where the distribution of the underlying random variable is discrete and has finitely many atoms. They use the DC representation to study a financial risk-return portfolio selection problem with a $\Var$ constraint. Recent work in \cite{thormann2024boosted} further advanced DC function optimization for $\Var$ by using boosted algorithms to enhance solution accuracy in portfolio optimization contexts. 
Due to the non-smooth and non-convex nature of $\Var$, studies have also been done to estimate $\Var$ (see, for example, \cite{taylor2008estimating, hendricks1996evaluation, kuester2006value}), in this direction $\CVar$ has long served as the best and most popular convex overestimate for $\Var$, also allowing for the construction of a convex overestimate for the $\Var$ minimization problem \cite{kunzi2006computational, sarykalin2008value, rockafellar2000optimization}. 
In this paper, we investigate the use of integrated quantile functions to describe a novel, unifying risk measure. In particular, for a given loss distribution function $F$ and probability levels $\alpha, \gamma$ satisfying $\alpha < \gamma$, the proposed risk measure leverages integrated quantile functions to account for all losses $l$ such that $F(l) \in [\alpha, \gamma]$.

%\textcolor{red}{Meanwhile, integrated quantile functions have been investigated as tools to capture a continuous distribution of risk across quantiles, presenting a richer depiction of potential outcomes and mitigating some of the restrictions inherent to single-quantile measures \cite{gushchin2017integrated}} \textcolor{red}{For instance $\Var$ only account for the loss distribution at a given probability level $\alpha$, CVaR accounts for the losses exceeding $\Var_{\alpha}$. We are interested in the expectation of losses with F() $\in [\alpha,\gamma]$.}.

\textbf{Contribution of the paper:} This paper presents a novel two-stage optimization framework designed to improve risk measurement by incorporating integrated quantile functions \cite{gushchin2017integrated}. Assuming that the distribution of the underlying random variable $\BT(x)$ is discrete and has finitely many atoms, the introduced two-stage optimization framework inputs probability levels $\alpha, \gamma$ (satisfying $\alpha < \gamma$) and models a new risk measure (referred to as $\BE_{\alpha - \gamma}[\BT(x)]$), special cases of which are $\Var_{\alpha}\BT(x)$ and $\CVar_{\alpha}\BT(x)$ \cite{rockafellar2000optimization}. 
By leveraging cutting-edge non-linear optimization techniques, we aim to solve the NP-Hard problem of minimizing $\BE_{\alpha - \gamma}[\BT(x)]$. Consequently, we also provide a fresh perspective on minimizing other key risk metrics ($\Var_{\alpha}\BT(x)$ and $\CVar_{\alpha}\BT(x)$), by incorporating a broader range of quantile information. 
A significant focus of this study is on exploring the relationship between the proposed risk measure $\BE_{\alpha - \gamma}[\BT(x)]$ and the integrated quantile function of $\BT(x)$, as well as developing various convexification techniques within the framework for this NP-Hard problem. These techniques aid in the construction of lower and upper convex estimators for the NP-Hard problem and provide new, tighter bounds for the $\Var$ minimization problem. This, in turn, supports the development of safe approximations for chance-constrained optimization problems with $\Var$ constraints (see for reference \cite{jiang2024also, nemirovski2007convex, nemirovski2012safe}). To validate and demonstrate the practical effectiveness of our framework, we apply it to real-world network traffic engineering (NTE) problems. Computational comparisons on the minimization of $\Var$ in the context of NTE problems illustrate the robustness and precision of our proposed method in managing network uncertainties compared to existing estimation approaches. 

\textbf{Outline of the paper:} 
In Section \ref{section_proposed_risk_measure}, we describe the construction and interpretation of the proposed risk measure $\BE_{\alpha - \gamma}[\cdot]$, also referred to as the ``$\alpha - \gamma$ expectation.'' Section \ref{section_exp_minimization} discusses the NP-Hard problem of minimizing this risk measure and highlights that the Var and CVaR minimization problems are special cases. In Sections \ref{subsection_underestimate} and \ref{subsection_overestimate}, we utilize several advanced optimization techniques to construct both over- and under-estimates for the NP-Hard problem, as well as for the $\Var$ minimization problem, thereby improving upon existing estimation methods for $\Var$ minimization. Finally, in Section \ref{section_case_study}, we present numerical computations that demonstrate the effectiveness of the proposed approach in real-life traffic engineering problems.

\section{Proposed risk measure: Rationale and design}\label{section_proposed_risk_measure}

\begin{figure}[hbtp]
    \centering
    \includegraphics[scale = 0.6]{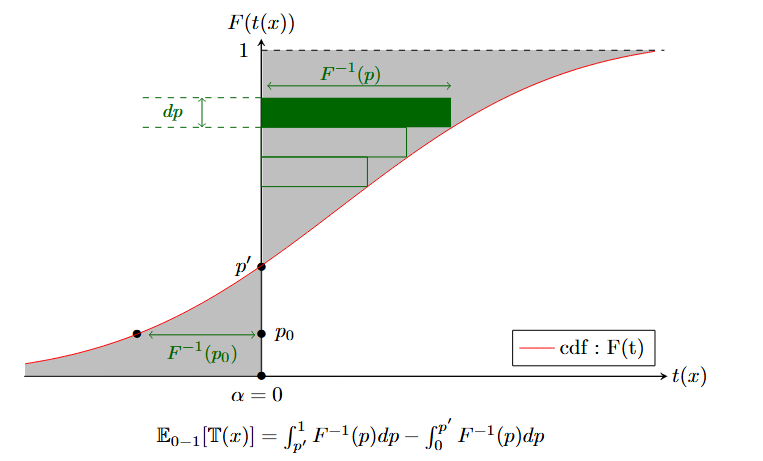}
    \caption{Computing $\BE_{0-1}[\BT(x)]$ via $F^{-1}(\cdot)$}
    \label{fig:exp_2}
\end{figure}

Given a nonempty set $\XP\subseteq \Re^{k}$, and a function $f:\XP\times \Re^{m} \rightarrow \Re$ such that for any $x\in \XP \subseteq \Re^{k}$, $\omega \in \Omega \subseteq \Re^{n}$ we have $f(x,\BY(\omega)) \in \Re$, where $\BY: (\Omega\subseteq \Re^{n},\mathcal{F}, \BP)  \rightarrow \Re^{m}$ is a random variable, $\mathcal{F}$ is a sigma algebra on $\Omega$, and $\BP$ is a probability measure defined on the space $(\Omega, \mathcal{F})$.
We assume that for any given realization of $\BY$, the function $f(x,\BY)$ is linear in $x$. For any given $x\in \XP$, we consider the random variable $f(x, \BY)$, and interchangeably refer to it as $\BT(x)$.
For the random variable $\BT(x)$, we let $F_{x}(\cdot):\Re \longrightarrow [0,1]$ denote its cumulative distribution function (cdf) and $F^{-1}_{x}(\cdot):[0,1]\longrightarrow \Re$ be the generalized inverse function of $F_{x}$.
For ease of notation, we will drop the subscript $x$ from $F_{x}(\cdot)$ and $F^{-1}_{x}(\cdot)$ when the context is clear.
For this setting, given a probability level $p$, the generalized inverse function $F^{-1}(\cdot)$ is defined as (\ref{eq:F_inverse}).
\begin{equation}\label{eq:F_inverse}
    F^{-1}(p) = \inf\bigl\{t : \BP\bigl(\BT(x) \leq t\bigr) = F(t) \geq p \bigr\}.  
\end{equation}
Given $x\in \XP$, and a probability level $\alpha$ satisfying $0 \leq \alpha < 1$, we leverage the definition of $F^{-1}(\cdot)$ in (\ref{eq:F_inverse}) to first define $\BE_{\alpha -1}[\BT(x)]$ as follows:
\begin{align}\label{eq:E_a_1_cont}
   \BE_{\alpha-1}[\BT(x)] =
    \begin{cases}
    \frac{1}{ (1-\alpha)} \int_{\alpha}^{1} F^{-1}(p) dp  & \text{ if } \alpha \geq p'
    \\
     \frac{1}{ (1-\alpha)}\Bigl[\int_{p'}^{1} F^{-1}(p) dp - \int_{\alpha}^{p'} F^{-1}(p) dp\Bigr] & \text{ if } \alpha < p',
    \end{cases}  
\end{align}
where $p'= \max_{p}\{p\in [0,1]: F^{-1}(p)\leq 0\}$. We also refer to the shaded region in Figure \ref{fig:exp_2}, that represents the quantity $\BE_{\alpha - 1}[\BT(x)]$, when $\alpha = 0$ and $F(\cdot)$ is the cdf of $\BT(x)$.
Given probability levels $\alpha, \gamma$ such that $\alpha < \gamma$, we take advantage of the definition of $\BE_{\alpha-1}[\BT(x)]$ in (\ref{eq:E_a_1_cont}) and define our proposed risk measure $\BE_{\alpha-\gamma}[\BT(x)]$ in  (\ref{eq:E_a_g_cont}). 
The proposed risk measure $\BE_{\alpha-\gamma}[\BT(x)]$, allows the upper probability level of $1$ in $\BE_{\alpha-1}[\BT(x)]$ to be replaced by $\gamma \leq 1$ such that $\alpha < \gamma$.
\begin{align}\label{eq:E_a_g_cont}
   \BE_{\alpha-\gamma}[\BT(x)] =
    \begin{cases}
    \frac{1}{ (\gamma -\alpha)} \int_{\alpha}^{\gamma} F^{-1}(p) dp  & \text{ if } \alpha \geq p'
    \vspace{0.3em}
    \\
     \frac{1}{ (\gamma -\alpha)}\Bigl[\int_{p'}^{\gamma} F^{-1}(p) dp - \int_{\alpha}^{p'} F^{-1}(p) dp\Bigr] & \text{ if } \alpha < p' \leq \gamma
     \vspace{0.3em}
     \\
     \frac{1}{ (\gamma -\alpha)}\Bigl[ - \int_{\alpha}^{\gamma} F^{-1}(p) dp\Bigr] & \text{ if } \alpha < \gamma < p',
    \end{cases}   
\end{align}
where $p'$ is as in (\ref{eq:E_a_1_cont}). We will hereafter interchangeably refer to $\BE_{\alpha-\gamma}[\BT(x)]$ as $\alpha-\gamma$ expectation of $\BT(x)$, and we will also make the following assumption about $\alpha, \gamma$.
\begin{enumerate}[label=\textbf{A1},align=left]
  \item \label{assumption:alpha_gamma} Probability levels $\alpha, \gamma$ satisfy the relation $\alpha < \gamma$.
\end{enumerate}

\subsection{Interpreting $\alpha - \gamma$ expectation: $\BE_{\alpha-\gamma}[\cdot]$}\label{subsection_Methods for expectation calculation}

In this section, given probability levels $\alpha, \gamma$ that satisfy the assumption \ref{assumption:alpha_gamma}, we further explore the interpretation of the risk measure $\BE_{\alpha - \gamma}[\BT(x)]$ in connection with integrated quantile functions. We also propose a two-stage max-min optimization formulation to model $\BE_{\alpha - \gamma}[\BT(x)]$.

\begin{definition}\label{def_IQF}
\textit{Integrated quantile function\cite{gushchin2017integrated}}
\item For a given $x\in \XP$, consider a random variable $\BT(x)$ with cdf $F(\cdot)$, the integrated distribution function of $\BT(x)$ is defined as $J(t) = \int_{0}^{t} F(x)dx$, and the \textit{integrated quantile function} (IQF) of $\BT(x)$ is defined as $K(p) = \sup_{t\in \Re}\{t p - J(t) \}$ for $p\in \Re$. 
Moreover, for any probability level $p'\in [0,1]$, the IQF of $\BT(x)$ at $p'$ \ie{} $K(p')$ is obtainable as $\int_{\widetilde{p}}^{p'}q(s)ds$, where $q(\cdot)$ is any quantile function (see \cite{gushchin2017integrated, embrechts2013note} for reference to quantile functions) of $\BT(x)$ and $\widetilde{p}$ is any zero of $\sup_{t\in \Re}\{t p - J(t) \}$.
\end{definition}

From \cite{gushchin2017integrated, embrechts2013note} we note that the generalized inverse function of $F(\cdot)$ as described in (\ref{eq:F_inverse}) is a quantile function of $\BT(x)$, and $\widetilde{p} = 0$ is a solution for $\sup_{t\in \Re}\{t p - J(t)\}$. Thus, we can equivalently write $K(p')$ in Definition \ref{def_IQF} as $\int_{0}^{p'}F^{-1}(s)ds$.
For any given probability level $p\in [0,1]$, hereafter we will use $\IQF(p) \mathord{:=} \int_{0}^{p}F^{-1}(s)ds$ to denote the IQF of $\BT(x)$ at $p$. 
Figure \ref{fig:IQF} shows the graph of a typical IQF of a random variable as a function of $p \in [0,1]$. Given any probability levels $\alpha, \gamma$, satisfying Assumption \ref{assumption:alpha_gamma}, $\BE_{\alpha-\gamma}[\BT(x)]$ is obtained as the slope of the line segment joining $(\alpha, \IQF(\alpha))$ and $(\gamma, \IQF(\gamma))$. 
This relation follows from our definition of $\BE_{\alpha - \gamma}[\BT(x)]$ in (\ref{eq:E_a_g_cont}) and noting that slope of the line segment joining $(\alpha, \IQF(\alpha))$ and $(\gamma, \IQF(\gamma))$ is obtained as:
\begin{align*}
\frac{1}{\gamma - \alpha}\bigl(\IQF(\gamma) - \IQF(\alpha)\bigr) 
& = \frac{1}{\gamma - \alpha}\Bigl(\int_{0}^{\gamma}F^{-1}(s)ds - \int_{0}^{\alpha}F^{-1}(s)ds  \Bigr)
\\
& = \frac{1}{\gamma - \alpha}\int_{\alpha}^{\gamma}F^{-1}(s)ds.
\end{align*}
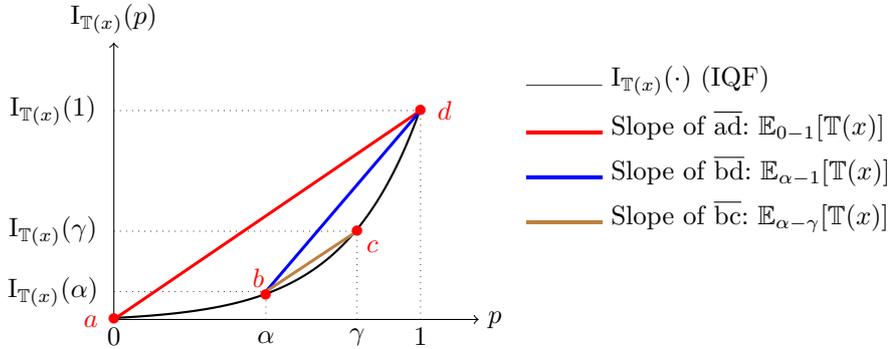
\begin{figure}[hbtp]
\centering
\begin{tikzpicture}[scale = 1, domain=0:4]
  % Axes
  \draw[->] (0,0.03) -- (4.8,0.03) node[right] {$p$};
  \draw[->] (0,0.03) -- (0,3.7) node[above] {$\IQF(p)$};

  % IQF curve
  \draw[domain=0:4.05, samples=100, smooth, variable=\x, black, thick]plot (\x, {0.05*exp(\x)});

  % Dotted lines for key points
  \draw[dotted, line width = 0.1mm] (4.035,0.03) -- (4.035,2.8);
  \draw[dotted, line width=0.1mm] (2,0.03) -- (2,0.4);
  \draw[dotted, line width=0.1mm] (3.2,0.03) -- (3.2,1.2);
  \draw[dotted, line width=0.1mm, black] (0.0, 0.4) -- (2, 0.4);
  \draw[dotted, line width=0.1mm, black] (0.0, 1.2) -- (3.2, 1.2);
  \draw[dotted, line width=0.1mm, black] (0.0, 2.8) -- (4.035, 2.8);

  % Slopes and line segments
  \draw[-, line width=0.4mm, red] (0,0.03) -- (4.035,2.8);
  \draw[-, line width=0.4mm, blue] (2.0, 0.4) -- (4.035,2.8);
  \draw[-, line width=0.4mm, brown] (2.0, 0.4) -- (3.2, 1.2);

  % Key points and annotations
  \node at (4.035, -0.2){$1$};
  \node at (0, -0.2){$0$};
  \node at (2.0, -0.2){$\alpha$};
  \node at (3.2, -0.2){$\gamma$};
  
  \node at (-0.8, 2.8){\textcolor{black}{$\IQF(1)$}};
  \node at (-0.8, 1.2){\textcolor{black}{$\IQF(\gamma)$}};
  \node at (-0.8, 0.4){\textcolor{black}{$\IQF(\alpha)$}};

  \node at (2,0.35){\textcolor{red}{$\bullet$}};
  \node at (0,0.03){\textcolor{red}{$\bullet$}};
  \node at (3.2, 1.2){\textcolor{red}{$\bullet$}};
  \node at (4.035,2.8){\textcolor{red}{$\bullet$}};
  
  \node at (-0.3,0){\textcolor{red}{$a$}};
  \node at (1.9,0.58){\textcolor{red}{$b$}};
  \node at (3.4, 1){\textcolor{red}{$c$}};
  \node at (4.35,2.8){\textcolor{red}{$d$}};

  % Legend outside the plot area (to the right) with spacing and thicker lines
  \node[anchor=west] at (5.3, 3.2) {\textcolor{black!85}{\textbf{\rule{1cm}{0.1mm}}} $\IQF(\cdot)$ (IQF)};
  \node[anchor=west] at (5.3, 2.6) {\textcolor{red}{\textbf{\rule{1cm}{0.5mm}}} Slope of $\overline{\rm ad}$: $\BE_{0 - 1}[\BT(x)]$};
  \node[anchor=west] at (5.3, 2.0) {\textcolor{blue}{\textbf{\rule{1cm}{0.5mm}}} Slope of $\overline{\rm bd}$: $\BE_{\alpha - 1}[\BT(x)]$};
  \node[anchor=west] at (5.3, 1.4) {\textcolor{brown}{\textbf{\rule{1cm}{0.5mm}}} Slope of $\overline{\rm bc}$: $\BE_{\alpha - \gamma}[\BT(x)]$};

\end{tikzpicture}

\caption{Interpreting $\BE_{p-p'}[\BT(x)]$ via Integrated Quantile Function of $\BT(x)$}

\label{fig:IQF}
\begin{tablenotes}
       \scriptsize 
       \footnotesize
       \item \textbf{[1]} Slope $\overline{\rm bd}$: Slope of the line joining the points $b = (\alpha, \IQF(\alpha))$ and $d = (1, \IQF(1))$. 
\end{tablenotes}
\end{figure}
We further note that for any probability level $\widetilde{p}$,  $\Var_{\widetilde{p}}\BT(x)$ as defined in (\ref{eq:traditional_var_def}) (see, for reference, \cite{nemirovski2007convex, rockafellar2000optimization, wozabal2012value}) is equivalent to the slope of the tangent drawn to the function $\IQF(\cdot)$ at the point $(\widetilde{p}, \IQF(\widetilde{p}))$.
\begin{equation}\label{eq:traditional_var_def}
\Var_{\widetilde{p}}\BT(x) = \inf_{\eta \in \Re}\{\eta: \BP(\BT(x) \leq \eta) = F(\eta) \geq \alpha\}.
\end{equation}
Hereafter, for ease of representation, we use the following notations.
\\
\textbf{Notation:} (i) Given the IQF of $\BT(x)$ \ie{} $\IQF(\cdot)$, and probability levels $p, p'$ satisfying $p<p'$, we will refer to $\text{S}(p - p')$ as the slope of the line segment joining $(p, \IQF(p))$ and $(p', \IQF(p'))$. We will also refer to $\text{S}(p)$ as the slope of the tangent to $\IQF(\cdot)$ at $(p, \IQF(p))$. (ii) $[n]$ would refer to the set of natural numbers $\{1, 2\dots,n\}$.

\begin{figure}[hbtp]
\centering
\begin{tikzpicture}[scale = 0.96, domain=0:4]
    \draw[->] (-4, 0) -- (6.8, 0) node[right] {$t(x)$};
    \draw[->] (0, 0) -- (0, 5.64) node[above] {$F(t(x))$};
    
    \node at (0.3,0.49){\textcolor{red}{$\alpha$}};
    
    \node at (-0.3,4.3){\textcolor{cyan}{$\gamma$}};

    \path[fill=pink] (-3, 0.49) --(-0.009, 0.49) -- (-0.009,1) -- (-3, 1) -- cycle;
    \path[fill=pink] (-0.945,1) -- (-0.009,1) -- (-0.009,2) -- (-0.945, 2) -- cycle;
    \path[fill=pink] (0.006, 2.01) -- (2.08,2.01) -- (2.08,3.0) -- (0.006, 3.0) -- cycle;
    \path[fill=pink] (0, 4.3) -- (4.07,4.3) -- (4.07,3) -- (0, 3) -- cycle;

    \draw[dashed, line width=0.35mm, blue] (-0.945,1) -- (0,1);
    \draw[dashed, line width=0.3mm, blue] (-3,0.49) -- (0,0.49);
    \draw[dashed, line width=0.35mm, blue] (2.08,2) -- (2.08,0);
    \draw[dashed, line width=0.35mm] (0,5) -- (4.0,5);
    \draw[dashed, line width=0.35mm, blue] (4.08,3) -- (4.08,0);
    \draw[dashed, line width=0.35mm, blue] (0, 3) -- (2.8, 3);
    \draw[dashed, line width=0.3mm, blue] (0,4.3) -- (4.07,4.3);

    \draw[-, line width=0.25mm, red] (-4, 0) -- (-3,0);
    \draw[-, line width=0.25mm, red] (-3, 1) -- (-1,1);
    \draw[-, line width=0.25mm, red] (2.08, 3) -- (4,3);
    \draw[-, line width=0.25mm, red] (4.08, 5) -- (6.5, 5);
    \draw[-, line width=0.25mm, red] (-1, 2) -- (2,2);

    \node at (-0.95,1){\textcolor{red}{$\circ$}};
    \node at (-0.95,2){\textcolor{red}{$\bullet$}};
    \node at (-3,0){\textcolor{red}{$\circ$}};
    \node at (-3,1){\textcolor{red}{$\bullet$}};
    \node at (4.08,3){\textcolor{red}{$\circ$}};
    \node at (4.08,5){\textcolor{red}{$\bullet$}};
    %\node at (-0.3, 5) {$1$};
    \node at (2.08,2){\textcolor{red}{$\circ$}};
    \node at (2.08,3){\textcolor{red}{$\bullet$}};
    
    \node at (0.3, 1){$p_1$};
    \node at (-3, -0.3){$t_1$};
    
    \node at (-1, -0.3){$t_2$};
    \node at (-0.3, 2.2){$p_2$};
    
    %\node at (0.6, 0.3){$p_0 = 0$};

    \node at (2.1, -0.3){$t_3$};
    \node at (-0.3, 3){$p_3$};
    
    \node at (4.1, -0.3){$t_4$};
    \node at (-0.3, 5){$1$};
    \node at (5.5, 0.5) {{\boxed{\textcolor{red}{\textbf{-----}} \text{ cdf : }$F(t)$}}};
\end{tikzpicture}
\caption{Shaded region: $(\gamma - \alpha)\BE_{\alpha-\gamma}[\BT(x)]$}
\label{fig:exp5}
\end{figure}
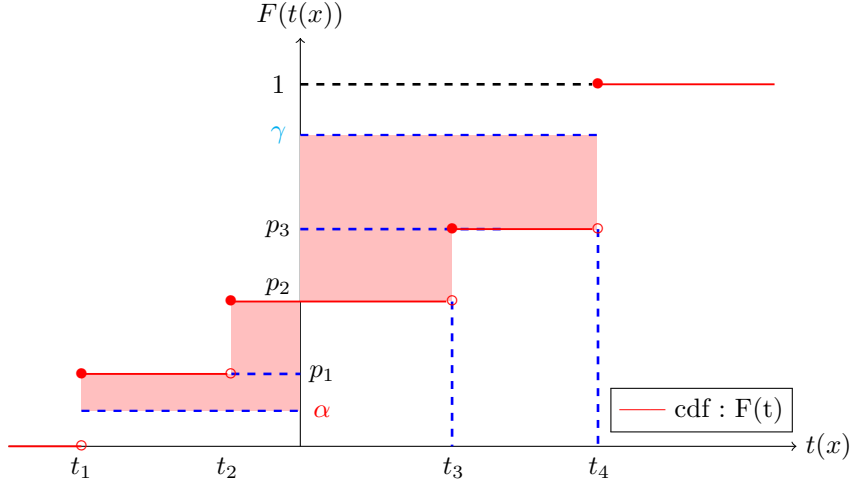

%From the above discussion, we note that the notion of the risk measure $\BE_{\alpha - \gamma}[\BT(x)]$ is valid for general random variables $\BT(x)$, however, we will hereafter restrict our analysis to the case when $\BT(x)$ for all $x\in \XP$ is a discrete random variable. To this, we will make the following assumption.
%
For our framework, in this paper, we will restrict our analysis to the case when $\BT(x)$ for all $x\in \XP$ is a discrete random variable. To this, we will make the following assumption.
\begin{enumerate}[label=\textbf{A2},align=left]
\item \label{assumption_discrete} $\BY: \Omega \rightarrow \Re^{m}$ as defined earlier is a discrete random vector having $n$ distinct realizations $y_1, \dots, y_n$ each in $\Re^{m}$ with the corresponding probabilities $p_1, \dots, p_n$.
\end{enumerate}
For any given $x \in \XP$, and for any realization $y_i$ for $i\in [n]$, we represent $f(x, y_i)$ as $T_{i}(x)$. Here, given $x\in \XP$, $\{T_{i}(x)\}_{i=1}^{n}$ is the set of realizations for the random variable $\BT(x) = f(x,\BY)$ corresponding to the $n$ realizations $\{y_i\}_{i=1}^{n}$ of $\BY$. Following Assumption \ref{assumption_discrete}, for all $i\in [n]$, $T_{i}(x)$ has an associated probability of $p_i$ that is induced by $y_i$. 
Now, for any predefined probability level $\alpha<1$, we compute $\BE_{\alpha-1}[\BT(x)]$ using the following linear program:
\begin{align}\label{LP_for_E_a1}
     \BE_{\alpha-1}\BT(x) = & \frac{1}{(1-\alpha)} \max_{w}\Bigl\{ \sum_{i\in [n]} T_{i}(x) w_i p_i \bigm| w \in \WP\Bigr\}, 
     \\
     & \text{where } \WP = \Bigl\{w \in [0,1]^{n} : \sum_{i\in [n]} w_i p_i = 1 - \alpha\Bigr\}. \nonumber
\end{align}
\begin{align}\label{LP_for_E_g1}
-\BE_{\gamma-1}\BT(x) = & \frac{1}{(1-\gamma)}\min_{w'}\Bigl\{-\sum_{i\in [n]} T_{i}(x) w'_i p_i: w'\in \WprimeP\Bigr\}
, 
\\
& \text{ where } \WprimeP = \Bigl\{w' \in [0,1]^{n} : \sum_{i\in [n]} w'_i p_i = 1 - \gamma\Bigr\}.\nonumber
\end{align}
Leveraging the quantities $\BE_{\alpha-1}[\BT(x)]$ and $\BE_{\gamma-1}[\BT(x)]$, as defined in (\ref{LP_for_E_a1}) and (\ref{LP_for_E_g1}) respectively, we then compute $\BE_{\alpha-\gamma}[\BT(x)]$ as in (\ref{eq:Exp_comp_alt_method}). Following Assumptions \ref{assumption:alpha_gamma} and \ref{assumption_discrete}, the shaded region in Figure \ref{fig:exp5}, describes the quantity $(\gamma - \alpha)\BE_{\alpha - \gamma}[\BT(x)]$ for one such typical instance. 
\begin{equation}\label{eq:Exp_comp_alt_method}
    \BE_{\alpha-\gamma} [\BT(x)] = \frac{1}{\gamma - \alpha} \Bigl((1-\alpha) \BE_{\alpha-1}[\BT(x)] - (1-\gamma) \BE_{\gamma - 1}[\BT(x)]\Bigr).   
\end{equation}
%
%\begin{figure}[!h]
%\centering
%\includegraphics[width=22em]{fig5.PNG}
%\caption{Shaded region: $(\gamma - \alpha)\BE_{\alpha-\gamma}[\BT(x)]$}
%\label{fig:exp5}
%\end{figure}
%
%
For a given $x\in \XP$, and known realizations for $\BT(x) = f(x,\BY)$ \ie{} $\{T_i(x)\}_{i=1}^{n}$, we use (\ref{LP_for_E_a1}) and (\ref{LP_for_E_g1}) to equivalently rewrite $\BE_{\alpha- \gamma}[\BT(x)]$ in (\ref{eq:Exp_comp_alt_method}) as the following two stage $\max$-$\min$ problem, %
%
%
%using \textcolor{red}{the $n$ realizations of $\BT(x) = f(x,\BY)$ \ie{} $\{T_i(x)\}_{i=1}^{n}$, equations} (\ref{LP_for_E_a1}) and (\ref{LP_for_E_g1}), we can equivalently rewrite $\BE_{\alpha- \gamma}[\BT(x)]$ in (\ref{eq:Exp_comp_alt_method}) as the following two stage $\max$-$\min$ problem, %
%
\begin{align}\label{eq:Exp_as_max_min_prob}
    \Gamma: \BE_{\alpha- \gamma}[\BT(x)] =  \max_{w\in \WP} \min_{w' \in \WprimeP} g(w,w'), \text{ where } g(w,w') \text{ is as follows,}
    \\
    g(w,w') = \frac{1}{(\gamma-\alpha)} \Bigl(\sum_{i\in [n]} T_i(x) w_i p_i - \sum_{i\in [n]} T_i(x) w'_i p_i\Bigr).\nonumber   
\end{align}
We note that $\WP$ and $\WprimeP$ as defined in (\ref{LP_for_E_a1}) and (\ref{LP_for_E_g1}) respectively are convex and compact sets. Further, $g: \WP \times \WprimeP \rightarrow \Re$ is a continuous function, such that for any fixed $w'\in \WprimeP$ and $w\in \WP$, $g(\cdot,w'):\WP \rightarrow \Re$ and $g(w, \cdot):\WprimeP \rightarrow \Re$ are respectively linear in $w$ and $w'$. 
Thus, $\max_{w\in \WP} \min_{w'\in \WprimeP} g(w,w') = \min_{w'\in \WprimeP} \max_{w\in \WP} g(w,w')$, which follows from von Neumann's minimax theorem (see, for instance, \cite{joo1980simple}).
We will therefore equivalently compute $\BE_{\alpha- \gamma}[\BT(x)]$ by solving $\min_{w'\in \WprimeP} \max_{w\in \WP} g(w,w')$.
Having already fixed an $x\in \XP$, for a given $w'^* \in \WprimeP$, we now look at the inner maximization problem in $\min_{w'\in \WprimeP} \max_{w\in \WP} g(w,w'^*)$, and let $\mathcal{D}(x,w'^*)$ represent its dual formulation. 
Then $\mathcal{D}(x,w'^*)$ is obtained as in (\ref{eq:dual_inner_max}). Here $(x,w'^*)$ in the notation of $\mathcal{D}(x,w'^*)$ indicates that the dual is obtained for a given $x, w'^*$.
\begin{equation}\label{eq:dual_inner_max}
\mathcal{D}(x,w'^*) \coloneqq \frac{1}{(\gamma-\alpha)}\min_{s, \theta} \Bigl\{(1-\alpha)s + \sum_{i\in [n]} (\theta_i - T_i(x) p_i w'^*_i)
\bigm| (s, \theta) \in \STpolytope(x,T(x)) \Bigr\},
\end{equation}
where $T(x) = (T_1(x),\dots,T_n(x))$ and $\STpolytope(x,T(x))$ is obtained as in (\ref{eq:STpolytope}).
We note that $x, T(x)$ appearing in the notation of $\STpolytope(x,T(x))$ in (\ref{eq:dual_inner_max}) is indicative of the fact that the polytope $\STpolytope(x,T(x))$ is obtained for a given $x\in \XP$, and known realizations $\{T_i(x)\}_{i=1}^{n}$ of $\BT(x)$.
\begin{equation}\label{eq:STpolytope}
    \STpolytope(x,T(x)) = \Bigl\{(s,\theta) \bigm| \theta_i + sp_i \geq p_i T_i(x), \ \theta_i \geq 0 \ \forall i \in [n]\Bigr\}.    
\end{equation}
We further note that, for any given $x\in \XP$ and $w' \in \WprimeP$, the choice of $s = \max_i T_i(x)$, and $\theta_i = 0$ for all $i\in [n]$ forms a feasible solution to $\mathcal{D}(x,w'^*)$.
Thus, the optimal value of $\mathcal{D}(x,w'^*) = \max_{w\in \WP} g(w,w'^*)$ which follows due to zero linear programming (LP) duality gap.
Now, from $\Gamma$ as defined in (\ref{eq:Exp_as_max_min_prob}), and leveraging the above LP dual construction, we can equivalently rewrite $\Gamma: \BE_{\alpha - \gamma}[\BT(x)] = \min_{w'\in \WprimeP} \mathcal{D}(x,w')$ as,
\begin{align}\label{eq:dual_E_a_g}
  \BE_{\alpha - \gamma}[\BT(x)] = \min_{s, \theta, w'} \quad & \frac{1}{(\gamma-\alpha)} \Bigl((1-\alpha)s + \sum_{i\in [n]} \theta_i - \sum_{i\in [n]} T_i(x) p_i w'_i\Bigr)
  \\
  \text{s.t. }  & (s, \theta) \in \STpolytope(x,T(x)), w' \in \WprimeP, \nonumber  
\end{align}
where $\STpolytope(x,T(x))$ and $\WprimeP$ are as defined in (\ref{eq:STpolytope}) and (\ref{LP_for_E_g1}) respectively. In the subsequent section, we explore the properties of the $\alpha-\gamma$ expectation of $\BT(x)$.

\subsection{Properties of $\alpha-\gamma$ expectation}
\begin{proposition}\label{proposition_Min_var_relation}
Given $x\in \XP$ and probability levels $\alpha, \gamma$ satisfying Assumption \ref{assumption:alpha_gamma}, for the random variable $\BT(x)$ we have $\displaystyle \Var_{\gamma}\BT(x) = \lim_{\alpha \uparrow {\gamma}} \BE_{\alpha - \gamma}[\BT(x)]$.
\end{proposition}
\proof
To prove the result, we first show that for any given probability levels $\alpha$ and $\gamma$ satisfying $\alpha<\gamma$, we have $\Var_{\alpha}\BT(x) \leq \BE_{\alpha-\gamma}[\BT(x)] \leq \Var_{\gamma}\BT(x)$ \ie{}
\begin{enumerate}
\item  $\Var_{p}\BT(x)$ is monotonically non-decreasing as a function of $p$, 
\item  $\BE_{\alpha-\gamma}[\BT(x)] \in [\Var_{\alpha}\BT(x), \Var_{\gamma}\BT(x)]$.
%$\Var_{\alpha}\BT(x) \leq \BE_{\alpha-\gamma}[\BT(x)]$ and $\Var_{\gamma}\BT(x) \geq \BE_{\alpha-\gamma}[\BT(x)]$.
\end{enumerate}
%in $p\in (0,1)$. 
%
The relation $\Var_{\alpha}\BT(x) \leq \Var_{\gamma}\BT(x)$ follows from the definition of $\Var_{*}\BT(x)$ in (\ref{eq:traditional_var_def}). 
To see this let $\eta_0$ be such that $\BP(\BT(x) \leq \eta_0) \geq \gamma$, since $\gamma > \alpha$ so $\eta_0$ also satisfies $\BP(\BT(x) \leq \eta_0) \geq \alpha$. 
Thus, we have $\{\eta_0 : \BP(\BT(x) \leq \eta_0) \geq \gamma\} \subseteq \{ \eta_0: \BP(\BT(x) \leq \eta_0) \geq \alpha\}$ and hence $\min_{\eta_0}\{\eta_0 : \BP(\BT(x) \leq \eta_0) \geq \gamma\} \geq \min_{\eta_0}\{ \eta_0: \BP(\BT(x) \leq \eta_0) \geq \alpha\}$, and 
$\Var_{\alpha}\BT(x) \leq \Var_{\gamma}\BT(x)$ follows.
We will now show that $\BE_{\alpha-\gamma}[\BT(x)]\in [\Var_{\alpha}\BT(x), \Var_{\gamma}\BT(x)]$ or equivalently $\BE_{\alpha-\gamma}[\BT(x)]\in [F^{-1}(\alpha), F^{-1}(\gamma)]$ which is true because of the equivalence of the definitions of $\Var_{*}\BT(x)$ and $F^{-1}(*)$ as in (\ref{eq:traditional_var_def}) and (\ref{eq:F_inverse}) respectively.
%
%We described earlier in (\ref{eq:E_a_g_cont}) and (\ref{eq:E_a_g_discont}), the computation of $\BE_{\alpha-\gamma}[\BT(x)]$ when $\BT(x)$ is continuous and discrete random variable respectively.
%
%As described earlier for our settings, since $\BT(x)$ has a discrete cdf, $\BE_{\alpha-\gamma}[\BT(x)]$ is computed using (\ref{eq:E_a_g_discont}).
%
Now, as stated earlier, since $\Var_{p}\BT(x)$ is monotonically non-decreasing in $p$, so for all $p_i$ for $i\in \{0,\dots,l\}$ such that $\alpha \leq p_i \leq \gamma$, we have $F^{-1}(\alpha) \leq F^{-1}(p_i) \leq F^{-1}(\gamma)$.
Moreover, since $\sum_{i=1}^{l}(p_i - p_{i-1}) = \gamma - \alpha$ for $p_0 = \alpha$ and $p_l = \gamma$, so we have the following relation:
\begin{equation}\label{eq:Var_ag_exp1}
    \frac{1}{\gamma-\alpha}F^{-1}(\alpha)(\gamma - \alpha) 
    \leq 
    \BE_{\alpha-\gamma}[\BT(x)] 
    \leq 
    \frac{1}{\gamma-\alpha} F^{-1}(\gamma)(\gamma - \alpha),   
\end{equation}
where $\frac{1}{\gamma-\alpha}F^{-1}(\alpha)(\gamma - \alpha) = \Var_{\alpha}\BT(x)$, $\frac{1}{\gamma-\alpha}F^{-1}(\gamma)(\gamma - \alpha) = \Var_{\gamma}\BT(x)$, and $\frac{1}{\gamma - \alpha}\sum_{i=1}^{l}F^{-1}(p_i)(p_i - p_{i-1}) = \BE_{\alpha - \gamma}[\BT(x)]$ using (\ref{eq:E_a_g_cont}).
%
%
%$\BE_{\alpha-\gamma}[\BT(x)] = $ which is in [] since $F^{-1}(p_i) \forall i$ is upper bounded by $F^{-1}(\gamma)$ and lower bounded by $F^{-1}(\alpha)$ moreover $\sum_{p_i}(p_i - p_{i-1}) = \gamma - \alpha$.
%
%
%
%
%where $\BE_{\alpha-\gamma}[\BT(x)]$ is computed as in (\ref{eq:E_a_g_discont}).
%
Lemma 2.1 in \cite{Kampke2015} shows that for a cdf $F(\cdot)$ which is right continuous and non-decreasing, $F^{-1}(p)$ as defined in (\ref{eq:F_inverse}) is non-decreasing and left continuous for $p\in (0,1)$.
%
%Lemma 2.1 in \cite{Kampke2015} shows that if $F(\cdot)$ is a cumulative distribution function which is right continuous and non-decreasing then $F^{-1}(p)$ as defined in (\ref{eq:F_inverse}) is non decreasing and left continuous for all $p\in (0,1)$.
%
Thus for $\epsilon >0$, we have $\lim_{\epsilon \rightarrow 0} F^{-1}(p-\epsilon) = F(p)$ for all $p\in (0,1)$.
Now, the equivalence of $\Var_{p}(\cdot)$ and $F^{-1}(p)$ from (\ref{eq:traditional_var_def}) and (\ref{eq:F_inverse}) respectively, imply that given the random variable $\BT(x)$,
%, $\Var_{p}\BT(x)$ is left continuous for $p\in (0,1)$. Hence, 
we have, $\lim_{\epsilon \rightarrow 0} \Var_{\gamma-\epsilon}\BT(x) = \Var_{\gamma}\BT(x)$ and from (\ref{eq:Var_ag_exp1}) it follows that:
\begin{equation*}
    \lim_{\alpha \uparrow \gamma}\Var_{\alpha}\BT(x) = \Var_{\gamma}\BT(x)
    \leq 
    \lim_{\alpha \uparrow \gamma} \BE_{\alpha-\gamma}[\BT(x)] 
    \leq 
    \Var_{\gamma}\BT(x),
\end{equation*}
%
%\begin{equation}\label{eq:Var_ag_exp}
%    \lim_{\alpha \uparrow \gamma}\Var_{\alpha}\BT(x)
%    =
%    \Var_{\gamma}\BT(x)
    %
%    \leq 
%    \lim_{\alpha \uparrow \gamma} \BE_{\alpha-\gamma}[\BT(x)] 
    %
%    \leq 
%    \Var_{\gamma}\BT(x),
%\end{equation}
%
implying $\Var_{\gamma} \BT(x) = \lim_{\alpha \uparrow \gamma} \BE_{\alpha - \gamma}[\BT(x)]$. 
\endproof
\begin{remark}
The result in Proposition \ref{proposition_Min_var_relation} can also be understood using the integrated quantile functions described in Section \ref{subsection_Methods for expectation calculation}, and noting that the IQF, $\IQF(p)$ is convex in $p$ \cite{gushchin2017integrated}.
Figure \ref{fig_E_to_Var}, shows that $\Var_{p}\BT(x)$ \ie{} the slope of the tangent to $\IQF(\cdot)$ at $(p, \IQF(p))$ is monotonically non decreasing in $p$, and $\BE_{\alpha - \gamma}[\BT(x)] \in [\Var_{\alpha}\BT(x), \Var_{\gamma}\BT(x)]$ \ie{} $\text{S}(\alpha - \gamma) \in [\text{S}(\alpha), \text{S}(\gamma)]$.
%, where $\text{S}(\alpha)$ and $\text{S}(\gamma)$ are the slopes of the tangent at $(\alpha, \IQF(\alpha))$ and $(\gamma, \IQF(\gamma))$ respectively.
%slope of the line segment joining $(\alpha, \IQF(\alpha))$ and $(\gamma, \IQF(\gamma))$ lies in the interval $[\text{s}_1, \text{s}_2]$, where $\text{s}_1$ and $\text{s}_2$ are the slopes of the tangent at $(\alpha, \IQF(\alpha))$ and $(\gamma, \IQF(\gamma))$ respectively.
%
We further note that for probability levels $\alpha < \alpha_1 < \alpha_2\leq \gamma$, $\text{S}(\alpha - \gamma) < \text{S}(\alpha_1 - \gamma) < \text{S}(\alpha_2 - \gamma) \leq S(\gamma)$, and hence $\lim_{\alpha \uparrow \gamma} \BE_{\alpha - \gamma}[\BT(x)] = \Var_{\gamma}\BT(x)$ (refer to Figure \ref{fig_E_to_Var}). 
%slope of $\overline{\rm bd} <$ slope of $\overline{\rm b_1d} <$ slope of $\overline{\rm b_2d}\leq$ slope of $\overleftrightarrow{\rm t_d}$, and hence $\lim_{\alpha \uparrow \gamma} \BE_{\alpha - \gamma}[\BT(x)] = \Var_{\gamma}\BT(x)$ (refer to Figure \ref{fig_E_to_Var}).
\end{remark}
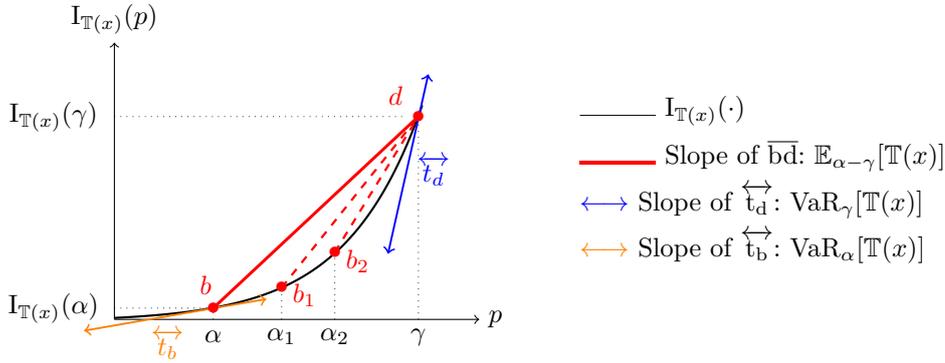
\begin{figure}[htbp]
\centering
\begin{tikzpicture}[scale=1, domain=0:4]
  % Axes
  \draw[->] (0,0.03) -- (4.8,0.03) node[right] {$p$};
  \draw[->] (0,0.03) -- (0,3.7) node[above] {$\IQF(p)$};

  % IQF curve
  \draw[domain=0:4.05, samples=100, smooth, variable=\x, black, thick]plot (\x, {0.05*exp(\x)});

  % Dotted lines for key points
  \draw[dotted, line width=0.1mm] (4, 0.03) -- (4,2.72);
  \draw[dotted, line width=0.1mm] (0, 2.72) -- (3.8,2.72);
  \draw[dotted, line width=0.1mm] (0, 0.18) -- (1.3,0.18);
  \draw[dotted, line width=0.1mm] (1.3, 0.03) -- (1.3, 0.18);
  \draw[dotted, line width=0.1mm] (2.2, 0.03) -- (2.2, 0.45);
  \draw[dotted, line width=0.1mm] (2.9, 0.03) -- (2.9, 0.92);

  % Lines and slopes
  \draw[-, line width=0.4mm, red] (1.3, 0.18) -- (4.0, 2.72);
  \draw[dashed, line width=0.3mm, red] (2.2, 0.45) -- (4, 2.72);
  \draw[dashed, line width=0.3mm, red] (2.9, 0.92) -- (4, 2.72);

  % Arrows representing variance
  \draw[<->, line width=0.25mm, blue] (3.6, 0.9) -- (4.13, 3.28);
  \node at (4.2,2){\textcolor{blue}{$\overleftrightarrow{t_d}$}};
  \draw[<->, line width=0.25mm, orange] (-0.4, -0.12) -- (2.0, 0.3);\node at (0.7,-0.3){\textcolor{orange}{$\overleftrightarrow{t_b}$}};

  % Key points and labels
  \node at (1.3, -0.2) {$\alpha$};
  \node at (2.2, -0.2) {$\alpha_1$};
  \node at (2.9, -0.2) {$\alpha_2$};
  \node at (4, -0.2) {$\gamma$};

  \node at (-0.8, 2.72) {$\IQF(\gamma)$};
  \node at (-0.8, 0.18) {$\IQF(\alpha)$};

  \node at (2.2, 0.45) {\textcolor{red}{$\bullet$}};
  \node at (2.9, 0.92) {\textcolor{red}{$\bullet$}};
  \node at (1.3, 0.18) {\textcolor{red}{$\bullet$}};
  \node at (4, 2.72) {\textcolor{red}{$\bullet$}};

  \node at (2.5, 0.35) {\textcolor{red}{$b_1$}};
  \node at (3.2, 0.8) {\textcolor{red}{$b_2$}};
  \node at (1.2, 0.48) {\textcolor{red}{$b$}};
  \node at (3.7, 3) {\textcolor{red}{$d$}};

  % Legend
  \node[anchor=west] at (6, 2.8) {\textbf{\rule{1cm}{0.1mm}} $\IQF(\cdot)$};
  \node[anchor=west] at (6, 2.2) {\textcolor{red}{\textbf{\rule{1cm}{0.4mm}}} Slope of $\overline{\rm bd}$: $\BE_{\alpha-\gamma}[\BT(x)]$};
  \node[anchor=west] at (6, 1.6) {\textcolor{blue}{$\longleftrightarrow$} Slope of $\overleftrightarrow{\rm t_d}$: $\Var_{\gamma}[\BT(x)]$};
  \node[anchor=west] at (6, 1.0) {\textcolor{orange}{$\longleftrightarrow$} Slope of $\overleftrightarrow{\rm t_b}$: $\Var_{\alpha}[\BT(x)]$};

\end{tikzpicture}

\caption{Proposition \ref{proposition_Min_var_relation} using Integrated Quantile Function for $\BT(x)$}

\label{fig_E_to_Var}
\begin{tablenotes}
       \scriptsize 
       \footnotesize
       \item \textbf{[1]} $\overleftrightarrow{t_b}, \overleftrightarrow{t_d}$: Tangent at the points $b = (\alpha, \IQF(\alpha))$ and $d = (\gamma, \IQF(\gamma))$ respectively. 
    \end{tablenotes}
\end{figure}
Following Assumption \ref{assumption_discrete}, we next propose Algorithm \ref{Algo_chose_alpha}, which given a probability level $\gamma$ identifies the threshold $\alpha^*$ such that for all $\alpha$, satisfying $\gamma > \alpha \geq \alpha^*$, we get $\BE_{\alpha - \gamma}[\BT(x)] = \Var_{\gamma}\BT(x)$.
\begin{algorithm}[htbp]
\caption{Compute $\alpha^*$ s.t. $\forall \alpha$ satisfying $\gamma > \alpha \geq \alpha^*$: $\BE_{\alpha-\gamma}[\BT(x) = \Var_{\gamma}\BT(x)$}\label{Algo_chose_alpha}
\textbf{Initialize:} $\alpha_k = 0$, $k=0$
\hspace*{6em} \textbf{Input:} $\gamma$, $\{p_i\}_{i=1}^{n}$, $\mathfrak{b}\in \BN$
\\
\quad \textbf{Output:} $\alpha^*$ such that for all $\alpha$ satisfying $\gamma > \alpha \geq \alpha^*$: $\BE_{\alpha-\gamma}[\BT(x)] = \Var_{\gamma}\BT(x)$ 
\begin{algorithmic}[1]
\Procedure{Finding-$\alpha^*$}{}
\\
    \label{step1-1}\hspace*{0.4cm} $ (O^*, \eta^*) \leftarrow \min_{\eta, O}\Bigl\{O \bigm| O \geq \sum_{i\in [n]}\eta_ip_i \geq \alpha_k, \eta \in \{0,1\}^{n}\Bigr\}$
    \\
    \label{step1-2}\hspace*{0.4cm} $k\leftarrow k+1$
    \While{$O^* < \gamma$}\label{step1-3}
    \State \label{step1-4} $\displaystyle\alpha_{k} \leftarrow \alpha_{k-1} + \frac{\gamma - \alpha_{k-1}}{\mathfrak{b}}$
    \State \label{step1-5}  $ (O^*, \eta^*) \leftarrow \min_{\eta, O}\Bigl\{O \bigm| O \geq \sum_{i\in [n]}\eta_ip_i \geq \alpha_k, \eta \in \{0,1\}^{n}\Bigr\}$
    %
    %\State \label{step1-6} $\eta^* \leftarrow \arg \min_{\eta}\biggl\{\sum_{i=1}^{n}\eta_ip_i \bigm| \alpha_k \leq \sum_{i=1}^{n}\eta_ip_i, \eta \in \{0,1\}^{n}   \biggr\}$
    %
    \State \label{step1-7} $k \leftarrow k+1$
    \EndWhile
    \State\Return $\alpha^* = \alpha_k$, $(O^*, \eta^*)$
\EndProcedure
\end{algorithmic}
\end{algorithm}
\begin{proposition}\label{prop_choose_epsilon} 
Given $x\in \XP$, let $T_1(x),\dots,T_n(x)$ represent the known $n$ realizations of the random variable $\BT(x)$, with associated probabilities $p_1,\dots,p_n$.
Given a probability level $\gamma$, if $\alpha$ is chosen such that $\gamma > \alpha \geq \alpha^*$, where $\alpha^*$ is obtained from Algorithm \ref{Algo_chose_alpha}, then $\BE_{\alpha-\gamma}[\BT(x)] = \Var_{\gamma}\BT(x)$.
\end{proposition}
\proof
Let $(O^*,\eta^*)$ and $\alpha^*$ be obtained at the termination of Algorithm \ref{Algo_chose_alpha}. 
%
%$(O^*,\eta^*)$ is obtained either at Steps \ref{step1-1} or \ref{step1-5}, where $O^*$ and $\eta^*$ are the optimal value and optimal solution of the optimization problem respectively.
%
%$\alpha^*$ is obtained using the update rule at Step \ref{step1-4} for some input $\mathfrak{b}\in \BN$.
%
Then, $(O^*,\eta^*)$ and $\alpha^*$ are such that $O^*,\eta^*$ are obtained as the optimal value and optimal solution of the optimization problem stated either at Step \ref{step1-1} or \ref{step1-5} of Algorithm \ref{Algo_chose_alpha}, and $\alpha^*$ is obtained using the update rule at Step \ref{step1-4} of Algorithm \ref{Algo_chose_alpha} for some input $\mathfrak{b}\in \BN$.
%\\
%\\
%From Algorithm \ref{Algo_chose_alpha} we note that $(O^*,\eta^*)$ and $\alpha^*$ obtained at the termination, are such that $(O^*,\eta^*)$ is obtained either at Steps \ref{step1-1} or \ref{step1-5}, and $\alpha^*$ is obtained using the update rule at Step \ref{step1-4} for some input $\mathfrak{b}\in \BN$.
%
Moreover at termination, the obtained $(O^*,\eta^*), \alpha^*$ satisfy the following relations:
\begin{equation}\label{eq:choosing_epsilon}
    \alpha^* \leq \gamma \text{ and } \alpha^* \leq O^* = \sum_{i\in [n]}\eta^{*}_{i} p_i \geq \gamma,    
\end{equation}
where $\alpha^* \leq \gamma$ follows from the update rule of $\{\alpha_k\}_{k\geq 1}$ in Step \ref{step1-4}, 
$\alpha^* \leq O^*$ follows from the constraint $\alpha_k \leq \sum_{i\in [n]} \eta_ip_i$ in the optimization problem solved at Steps \ref{step1-1} and \ref{step1-5},
and $O^* \geq \gamma$ follows from the termination criterion in Step \ref{step1-3}.
The relations in (\ref{eq:choosing_epsilon}) imply that $F^{-1}(p)$ or equivalently $\Var_{p}\BT(x)$ as a function of the probability level $p$, does not have any jumps for $p\in [\alpha^*, \gamma]$, thus $\Var_{\alpha^*}\BT(x) = \Var_{\gamma}\BT(x)$.
Now, for any $\alpha$ such that $\alpha^* \leq \alpha < \gamma$ we have:
\begin{equation}\label{eq:Exp_equiv_var}
    \Var_{\alpha^*}\BT(x) \leq \Var_{\alpha}\BT(x)\leq \BE_{\alpha-\gamma}[\BT(x)] \leq \Var_{\gamma}\BT(x),    
\end{equation}
where the first inequality follows as $\Var_{p}\BT(x)$ is non-decreasing in $p$, and the rest are already argued in (\ref{eq:Var_ag_exp1}) as part of the proof for Proposition \ref{proposition_Min_var_relation}.
Since, $\Var_{\alpha^*}\BT(x) = \Var_{\gamma}\BT(x)$ follows from (\ref{eq:choosing_epsilon}), hence equality holds throughout in (\ref{eq:Exp_equiv_var}).
%    $\Var_{\gamma}\BT(x) = \Var_{\alpha^*}\BT(x) \leq \Var_{\alpha}\BT(x)\leq \BE_{\alpha-\gamma}[\BT(x)] \leq \Var_{\gamma}\BT(x)$.
%
So for all $\alpha\geq \alpha^*$, when $\alpha < \gamma$, we have $\BE_{\alpha-\gamma}[\BT(x)] = \Var_{\gamma}\BT(x)$. 
\endproof
%
%
%
%\begin{align}\label{eq:knapsack_IP}
%    \epsilon_0 = \min_{\eta} \Bigl\{\sum_{i=1}^{n}\eta_i p_i - (1-\alpha) \bigm| \sum_{i=1}^{n} \eta_i p_i \geq 1 - \alpha,\ \eta \in \{0,1\}^{n}  \Bigr\},
%\end{align}
%
%
Given $x\in \XP$, in the subsequent parts of this section, we further explore properties of the $\alpha - \gamma$ expectation of $\BT(x)$.
\begin{lemma}\label{lemma_convex_combin_E}
Given $x\in \XP$, and probability levels $\alpha, \alpha_1, \gamma$ satisfying $\alpha < \alpha_1 < \gamma \leq 1$, the relation: $\BE_{\alpha-\alpha_1}[\BT(x)] < \BE_{\alpha-\gamma}[\BT(x)] < \BE_{\alpha_1-\gamma}[\BT(x)]$ holds true.
\end{lemma}
\proof
To prove the above result, we will show that $\BE_{\alpha-\gamma}[\BT(x)]$ is a convex combination of $\BE_{\alpha-\alpha_1}[\BT(x)]$ and $\BE_{\alpha_1-\gamma}[\BT(x)]$.
From (\ref{eq:Exp_comp_alt_method}) we note that $\BE_{\alpha-\alpha_1}[\BT(x)] = \frac{(1-\alpha) \BE_{\alpha-1}[\BT(x)] - (1-\alpha_1) \BE_{\alpha_1 - 1}[\BT(x)]}{\alpha_1 - \alpha}$, 
and
$\BE_{\alpha_1-\gamma}[\BT(x)] = \frac{(1-\alpha_1) \BE_{\alpha_1-1}[\BT(x)] - (1-\gamma) \BE_{\gamma - 1}[\BT(x)]}{\gamma - \alpha_1}$.
For $\lambda$, $1-\lambda$ chosen as $\frac{\alpha_1-\alpha}{\gamma-\alpha}>0$ and $\frac{\gamma-\alpha_1}{\gamma-\alpha}>0$ respectively, we have $\BE_{\alpha-\gamma}[\BT(x)] 
=
\lambda \BE_{\alpha-\alpha_1}[\BT(x)] + (1-\lambda)\BE_{\alpha_1-\gamma}[\BT(x)]$. 
%
%= \frac{(1-\alpha) \BE_{\alpha-1}[\BT(x)] - (1-\gamma) \BE_{\gamma - 1}[\BT(x)]}{\gamma - \alpha}$.
%
Thus, $\BE_{\alpha-\alpha_1}[\BT(x)] < \BE_{\alpha-\gamma}[\BT(x)] < \BE_{\alpha_1-\gamma}[\BT(x)]$ follows.
\endproof
\begin{corollary}\label{lemma_E_monotonic}
Given $x\in \XP$, for any probability levels $\alpha$ and $\alpha_1$ satisfying $\alpha < \alpha_1 < 1$,
%% $\gamma < 1$ because $\BE_{\gamma-1}\BT(x)$ is not defined if $\gamma = 1$, because denominator in the obj becomes 0.
the relation: $\BE_{\alpha -1}[\BT(x)] < \BE_{\alpha_1 - 1}[\BT(x)]$ holds true. 
\end{corollary}
\proof
The result follows by substituting $\gamma = 1$ in Lemma \ref{lemma_convex_combin_E}.
%We first recall the definition of $\BE_{\alpha-1}[\BT(x)]$ in (\ref{LP_for_E_a1}) and let $\bar{w}$ be any feasible solution to $\text{P}_{\alpha}:\frac{1}{1-\alpha} \max_{w}\Bigl\{\sum_{i\in [n]} T_i(x) w_i p_i \bigm| \sum_{i\in [n]} w_i p_i = 1 - \alpha, w \in [0,1]^{n}\Bigr\}$, resulting in an objective value of $\bar{o} = \frac{1}{1-\alpha}\sum_{i\in [n]}T_i(x)p_i\bar{w}_i$.
%
%Then, $\bar{v} = \frac{1-\gamma}{1-\alpha}\bar{w}$ is a feasible solution to 
%
%$\text{P}_{\gamma}: \frac{1}{1-\gamma}\max_{v}\Bigl\{\sum_{i\in [n]} T_i(x) v_i p_i \bigm| \sum_{i\in [n]} v_i p_i = 1 - \gamma, v \in [0,1]^{n}\Bigr\}$ also resulting in an objective value of $\bar{o}$.
%
%So any solution feasible to $\text{P}_{\alpha}$ can be scaled down by $\frac{1-\gamma}{1-\alpha}\in [0,1)$, to be feasible to $\text{P}_{\gamma}$. Thus, $\BE_{\alpha-1}[\BT(x)] \leq \BE_{\gamma-1}[\BT(x)]$ follows.
%being maximization problems $\BE_{\alpha-1}[\BT(x)] < \BE_{\gamma-1}[\BT(x)]$ follows.
\endproof
\begin{remark}
The results in Lemma \ref{lemma_convex_combin_E} and Corollary \ref{lemma_E_monotonic} can also be interpreted using the convexity of $\IQF(p)$ with respect to $p$. Figure \ref{fig_monotonicity_E} shows that for given probability levels $\alpha, \alpha_1, \gamma$ satisfying $\alpha < \alpha_1 < \gamma \leq 1$, we have $\text{S}(\alpha-\alpha_1) < \text{S}(\alpha-\gamma) < \text{S}(\alpha_1 - \gamma)$. Moreover, the result in Corollary \ref{lemma_E_monotonic} follows by choosing $\gamma = 1$ in Figure \ref{fig_monotonicity_E} \ie{} $\text{S}(\alpha - 1) < \text{S}(\alpha_1 - 1)$.
\begin{figure}[htbp]
\centering
\begin{tikzpicture}[scale=1, domain=0:4]
  % Axes
  \draw[->] (0,0.03) -- (4.8,0.03) node[right] {$p$};
  \draw[->] (0,0.03) -- (0,3.7) node[above] {$\IQF(p)$};

  % IQF curve
  \draw[domain= 0:4.02, samples=100, smooth, variable=\x, black, thick]plot (\x, {0.05*exp(\x)});

  % Dotted lines for key points
  \draw[dotted, line width=0.1mm] (4, 0.03) -- (4, 2.72);
  \draw[dotted, line width=0.1mm] (0, 2.72) -- (3.8, 2.72);
  \draw[dotted, line width=0.1mm] (0, 0.18) -- (1.3, 0.18);
  \draw[dotted, line width=0.1mm] (0, 0.92) -- (2.9, 0.92);
  \draw[dotted, line width=0.1mm] (1.3, 0.03) -- (1.3, 0.18);
  \draw[dotted, line width=0.1mm] (2.9, 0.03) -- (2.9, 0.92);

  % Lines and slopes
  \draw[-, line width=0.4mm, red] (1.3, 0.18) -- (4, 2.72);
  \draw[-, line width=0.3mm, blue] (2.9, 0.92) -- (4, 2.72);
  \draw[-, line width=0.3mm, orange] (1.3, 0.18) -- (2.9, 0.92);

  % Key points and labels
  \node at (1.3, -0.2) {$\alpha$};
  \node at (2.9, -0.2) {$\alpha_1$};
  \node at (4, -0.2) {$\gamma$};

  \node at (-0.8, 2.72) {$\IQF(\gamma)$};
  \node at (-0.8, 0.18) {$\IQF(\alpha)$};
  \node at (-0.8, 0.92) {$\IQF(\alpha_1)$};

  \node at (1.3, 0.18) {\textcolor{red}{$\bullet$}};
  \node at (2.9, 0.92) {\textcolor{red}{$\bullet$}};
  \node at (4, 2.72) {\textcolor{red}{$\bullet$}};

  \node at (1.2, 0.48) {\textcolor{red}{$b$}};
  \node at (3.2, 0.8) {\textcolor{red}{$b_2$}};
  \node at (3.7, 3) {\textcolor{red}{$d$}};

  % Legend
  \node[anchor=west] at (6, 2.8) {\textbf{\rule{1cm}{0.1mm}} $\IQF(\cdot)$};
  \node[anchor=west] at (6, 2.2) {\textcolor{red}{\textbf{\rule{1cm}{0.5mm}}} Slope of $\overline{\rm bd}$: $\BE_{\alpha - \gamma}[\BT(x)]$};
  \node[anchor=west] at (6, 1.6) {\textcolor{blue}{\textbf{\rule{1cm}{0.5mm}}} Slope of $\overline{\rm b_2d}$: $\BE_{\alpha_1 - \gamma}[\BT(x)]$};
  \node[anchor=west] at (6, 1.0) {\textcolor{orange}{\textbf{\rule{1cm}{0.5mm}}} Slope of $\overline{\rm bb_2}$: $\BE_{\alpha - \alpha_1}[\BT(x)]$};

\end{tikzpicture}

\caption{Lemma \ref{lemma_convex_combin_E} and Corollary \ref{lemma_E_monotonic} using $\IQF(\cdot)$}

\label{fig_monotonicity_E}
\end{figure}
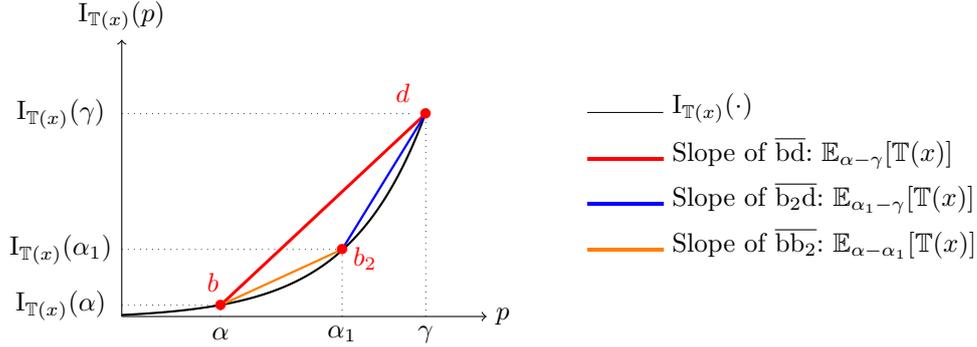
\end{remark}
\begin{lemma}\label{lemma_w_wp}
Given $x\in \XP$ and probability levels $\alpha,\gamma$ satisfying Assumption \ref{assumption:alpha_gamma}, consider the problem of computing $\BE_{\alpha-\gamma}[\BT(x)] = \max_{w\in \WP}\min_{w'\in \WprimeP} g(w,w')$ where $g(w,w') = \frac{\sum_{i\in [n]} T_i(x) p_i (w_i - w'_i)}{\gamma - \alpha}$ as described in (\ref{eq:Exp_as_max_min_prob}).
Any solution $(w^*\in \WP,w^{'*}\in \WprimeP)$ optimal to $\max_{w\in \WP}\min_{w'\in \WprimeP} g(w,w')$ satisfies $w_{i}^{*} \geq w_{i}^{'*}$ for all $i\in [n]$.
\end{lemma}
\proof
Given $x\in \XP$, realization of the random variable $\BT(x)$ \ie{} $\{T_i(x)\}_{i\in [n]}$, and its corresponding probabilities $\{p_i\}_{i\in[n]}$.
We assume w.l.o.g that $\{T_i(x)\}_{i\in [n]}$ can be sorted in a decreasing order to get a permutation $\sigma$ such that $T_{\sigma(1)}(x) \geq \dots \geq T_{\sigma(n)}(x)$.
If $w^{'*} = (w^{'*}_{\sigma(1)},\dots, w^{'*}_{\sigma(n)})$ is given to be an optimal solution to $\max_{w\in \WP}\min_{w'\in \WprimeP} g(w,w')$, then there exists $k\in [n]$ such that $w^{'*}_{\sigma(j)} = 1$ for all $j \leq k$,
$w^{'*}_{\sigma(j)} \in \{1,e\}$ for $j=k+1$ where $e\in (0,1)$, 
$w^{'*}_{\sigma(j)} = 0$ for all $j>k+1$, and $\sum_{i\in [k+1]}w^{'*}_{\sigma(i)} p_{\sigma(i)} = 1 - \gamma$.
Since $w^*$ is given to be optimal to $\max_{w\in \WP}\min_{w'\in \WprimeP} g(w,w')$, hence $w^* \in \arg\max_{w\in \WP}g(w,w^{'*})$.
We note that $w^*$ is obtained as $w^* = (w^{*}_{\sigma(1)}, \dots, w^{*}_{\sigma(k)}, w^{*}_{\sigma(k+1)}, w^{*}_{\sigma(k+2)} \dots, w^{*}_{\sigma(n)})$, where $w^{*}_{\sigma(j)} = w^{'*}_{\sigma(j)} = 1$ for all $j\leq k$, and $w^{*}_{\sigma(j)} \geq w^{'*}_{\sigma(j)}$ for all $j\geq k+1$. 
This follows since $\sum_{i\in [k]}w^{'*}_{\sigma(i)} p_{\sigma(i)} + w^{'*}_{\sigma(k+1)}p_{\sigma(k+1)} = \sum_{i\in [k]}w^{*}_{\sigma(i)} p_{\sigma(i)} + w^{'*}_{\sigma(k+1)}p_{\sigma(k+1)} = 1-\gamma < 1 - \alpha$.
On fixing $w' = w^{'*}$ in the objective in $\max_{w\in \WP} \min_{w'\in \WprimeP} g(w,w^{'*})$, which aims at maximizing $\frac{1}{\gamma - \alpha}\sum_{i\in [n]}(T_i(x)w_ip_i - T_i(x)w^{'*}_ip_i)$ as in (\ref{eq:Exp_as_max_min_prob}) over $w\in \WP$, forces $w_i$ to be picked corresponding to larger $T_i(x)$ values.
Thus $w^{*}_i \geq w^{'*}_i$ follows for all $i\in [n]$.
\endproof
Given $x\in \XP$, having described the proposed risk measure for the random variable $\BT(x)$, and discussed its properties, we are now interested in solving for an optimal $x, T(x) = (T_{i}(x))_{i=1}^{n}$, that minimizes $\BE_{\alpha-\gamma}[\BT(x)]$. In particular, we are interested in computing $\nu^*$ obtained as follows: 
\begin{equation}\label{eq:problem_definition}
    \nu^* = \min_{x, T(x)}\bigl\{\BE_{\alpha - \gamma}[\BT(x)] \mid x \in \XP, T_i(x) = f(x,y_i) \ \forall i \in [n]\bigr\}.    
\end{equation}
Computing $\nu^*$ is NP-hard. The non-convexity of computing $\nu^*$ follows from Propositions \ref{proposition_Min_var_relation} and \ref{prop_choose_epsilon}, and noting that the $\Var$-minimization problem: $\min_{x, T(x)} \bigl\{\Var_{\gamma} \BT(x)\mid x\in \XP, T_i(x) = f(x,y_i)\ \forall i \in [n] \bigr\}$ which is a special case of (\ref{eq:problem_definition}) is non-convex even when $\{(x,T(x)): x\in \XP, T_i(x) = f(x,y_i) \forall i \in [n]\}$ is a closed convex set (described for instance in \cite{larsen2002algorithms, pang2004global, wozabal2012value, kast1998var, sarykalin2008value}). Hence, we will be interested in estimating $\nu^*$.
In the subsequent sections, for concreteness in the discussion of computing $\nu^*$, we make the following assumption without loss of generality (w.l.o.g).
\begin{enumerate}[label=\textbf{A3},align=left]
\item \label{assumption_XT_polytope} Constraints defining $x\in \XP$ and $T_i(x) = f(x,y_i)$ for all $i\in [n]$ are represented as $Ax + BT(x) \leq c$ where $A\in \Re^{l\times k}, B\in \Re^{l\times n}$, $c\in \Re^{l}$, $x\in \Re^{k}$, and $T(x) = (T_i(x))_{i=1}^{n}$. In particular, we let $\XTP$ denote the $\{(x, T(x)) \mid Ax + BT(x) \leq c\}$, which is assumed to be a non-empty, closed and convex set.
\end{enumerate}

\section{Minimizing $\alpha - \gamma$ expectation}\label{section_exp_minimization}

%
%We next introduce approaches to compute $\nu^*$ as described in (\ref{eq:problem_definition}). 
%
We leverage the formulation for $\BE_{\alpha-\gamma}[\BT(x)]$ in (\ref{eq:dual_E_a_g}), and propose a bilinear optimization program in (\ref{eq:bilenear_formulation}) to compute $\nu^*$:
\begin{align}\label{eq:bilenear_formulation}
\nu^* = \min_{\psi}  \quad &  \frac{1}{\gamma- \alpha} \Bigl((1-\alpha)s + \sum_{i\in [n]}\theta_i - \sum_{i\in [n]}T_i(x) p_i w'_i\Bigr)
\\
\text{s.t. } &  (s,\theta) \in \STpolytope(x,T(x)),\  w' \in \WprimeP, \ (x,T(x))\in \XTP, \nonumber  
\end{align}
where $\STpolytope(x,T(x))$, $\WprimeP$, and $\XTP$ are as defined in (\ref{eq:STpolytope}), (\ref{LP_for_E_g1}), and Assumption \ref{assumption_XT_polytope}  respectively.
Moreover, the minimization in (\ref{eq:bilenear_formulation}) occurs over $\psi = (x,T(x),s,\theta,w')$. Next, in this section, we explore the properties of this bilinear optimization formulation.
%
%
%Hereafter we will also assume the following,
%
%\begin{assumption}\label{%assumption_XT_polytope}
%Constraints defining $x\in \XP$ and $T_i(x) = f(x,y_i)$ for all $i\in [n]$ are represented as $Ax + BT(x) \leq c$ where $A\in \Re^{l\times k}, B\in \Re^{l\times n}$ and $c\in \Re^{l}$. Moreover $\XTP: \{(x,T(x)) \mid Ax + BT(x) \leq c\}$ is non empty, convex and a bounded set.
%\end{assumption}
%
%The bilinear formulation computing $\nu^* = \min_{x\in \XP}\Gamma$ where $\Gamma$ is as defined in (\ref{eq:Exp_as_max_min_prob}) is obtained as,
%
%\begin{align}\label{eq:bilenear_formulation}
%\nu^* = \min_{\psi}  \quad &  \frac{(1-\alpha)s + \sum_{i\in [n]}\theta_i - \sum_{i\in [n]}T_i(x) p_i (1+u_i-z_i)}{\gamma- \alpha}
%\\
%& (s,\theta) \in \STpolytope(x,T(x)), \ w'=(u,z) \in \WprimeP \nonumber 
%\\  
%& (x,T(x))\in \XTP, \text{ for }\psi = (x,T(x),s,\theta,w')\nonumber 
%\end{align}
%
%
%
%As stated above, we compute $\nu^* = \min_{x\in \XP} \BE_{\alpha-\gamma}[\BT(x)]$, using the bilinear optimization problem in (\ref{eq:bilenear_formulation}). 
%
%In this section, we will further explore this problem.
%
\begin{proposition}\label{prop:CVar_a_special_case}
   Consider the problem defined in (\ref{eq:bilenear_formulation}). 
   Solving for $\nu^*$, when $w'_i = 0$ for all $i \in [n]$, is equivalent to solving $\min_{x,T(x)}\{\CVar_{\alpha} \BT(x)\mid (x,T(x)) \in \XTP$\}, where $\CVar_{\alpha}\BT(x)$ is defined as $\min_{\eta}\bigl\{\eta + \frac{1}{1-\alpha}\BE[\BT(x) - \eta]^{+}\bigr\}$ and $[t]^{+} = \max\{0,t\}$ (see for instance \cite{sarykalin2008value, rockafellar2000optimization}).
\end{proposition}

\proof
Consider the optimization problem in (\ref{eq:bilenear_formulation}) for a given $x\in \XP$, and known realizations $T(x) = (T_1(x),\dots,T_n(x))$ of $\BT(x)$. 
On fixing $w'_i = 0$ for all $i\in [n]$, we have $\gamma = 1$ from (\ref{LP_for_E_g1}), and the problem further reduces to $\Gamma(x)$ obtained as follows:
\begin{equation}\label{eq:gamma_w_0}
\begin{alignedat}{3}
    \Gamma(x) \mathord{:=}  \min_{s, \theta} \  & \frac{(1-\alpha)s + \sum_{i\in [n]} \theta_i}{1 - \alpha}
    \\
      & \theta_i + sp_i \geq T_i(x) p_i \ &  \forall i\in [n] 
      \\
      & \theta_i \geq 0 \ &  \forall i\in [n] 
\end{alignedat}
\quad 
\begin{alignedat}{3}
    \Gamma^{D}(x) \mathord{:=}  \max_{w} \  & \frac{\sum_{i\in [n]} T_i(x)p_i w_i}{1 - \alpha}
      \\
      & w \in [0,1]^{n} 
      \\
      & \sum_{i\in [n]} w_i p_i = 1 - \alpha 
\end{alignedat}
\end{equation}
Clearly, $\theta_i = 0$ for all $i\in [n]$ and $s= \max_{i} T_i(x)$ is a feasible solution to $\Gamma(x)$, and $\Gamma^{D}(x)$ being the dual of $\Gamma(x)$, we thus have $\Gamma^{*}(x) = \Gamma^{D^*}(x)$ where $\Gamma^{*}(x)$ and $\Gamma^{D^*}(x)$ represent the optimal value of $\Gamma(x)$ and $\Gamma^{D}(x)$ respectively. This follows due to the zero LP duality gap.
Further $\Gamma^{D^*}(x) = \BE_{\alpha - 1}\BT(x)$, which follows from (\ref{LP_for_E_a1}). Hence, fixing $w'_i = 0$ for all $i\in [n]$, $\nu^*$ in (\ref{eq:bilenear_formulation}) is equivalent to $\min_{x, T(x)\in \XTP} \BE_{\alpha -1}\BT(x)$. Now, we are left to show that $\BE_{\alpha -1}\BT(x) = \min_{\eta}\bigl\{\eta + \frac{1}{1-\alpha}\BE[\BT(x) - \eta]^{+}\bigr\}$. 
To do this we equivalently rewrite $\min_{\eta}\bigl\{\eta + \frac{1}{1-\alpha}\BE[\BT(x) - \eta]^{+}\bigr\}$ as $\mathfrak{P}$ below.
 \begin{equation}\label{eq:Cvar_recc_def}
\begin{alignedat}{3}
    \mathfrak{P} := \min_{\eta, \Psi} \ & \eta + \sum_{i\in [n]} \Psi_{i}
     \\
     & \Psi_i \geq \frac{\bigl(T_i(x) - \eta\bigr) p_i}{1-\alpha} \quad & \forall i \in [n]
     \\
     & \Psi_{i} \geq 0 \quad & \forall i \in [n]
\end{alignedat}
\quad
\begin{alignedat}{3}
     \mathscr{D} := \max_{d} \  & \sum_{i\in [n]} T_i(x)p_i d_i
      \\
      & d_i \in \Bigl[0,\frac{1}{1-\alpha}\Bigr] \quad &  \forall i\in [n] 
      \\
      & \sum_{i\in [n]} d_i p_i = 1.
\end{alignedat}
\end{equation}
The formulation $\mathscr{D}$ in (\ref{eq:Cvar_recc_def}) represents the dual formulation of $\mathfrak{P}$. 
We note that $\Psi_i = 0$ for all $i\in [n]$ and $\eta = \max_{i}T_i(x)$ is a feasible solution to $\mathfrak{P}$, so there is no duality gap. 
In $\mathscr{D}$, we then scale the variables $d_i$ for all $i\in [n]$, by multiplying its constraints with $(1-\alpha)$, and by defining $w_i \mathord{:=} d_i(1-\alpha)$. 
%
%We then scale the variables $D_i$ for all $i\in [n]$, in $\mathscr{D}$ by multiplying its constraints by $(1-\alpha)$, and defining $w_i \mathord{:=} D_i(1-\alpha)$. 
%
Following this we get,
\begin{equation*}
    \displaystyle \mathscr{D} :=  \max_{w}\biggl\{\frac{1}{1-\alpha}\sum_{i=1}^{n} T_i(x)p_i w_i \Bigm| \sum_{i=1}^{n} w_i p_i = 1 - \alpha,\ w_i \in [0,1] \quad \forall i \in [n] \biggr\}.
\end{equation*}
Finally, from (\ref{LP_for_E_a1}) we have $\BE_{\alpha -1}\BT(x) = \min_{\eta}\bigl\{\eta + \frac{1}{1-\alpha}\BE[\BT(x) - \eta]^{+}\bigr\}$, and hence the result follows.
\endproof
For a given probability level $\alpha$, we next examine that $\nu^*$ in (\ref{eq:bilenear_formulation}) is monotonic increasing for $\gamma \in (\alpha,1]$. 
%
%
%
%
\begin{comment}
\begin{corollary}\label{corollary_E_monotonic}
If probability levels $\beta$ and $\beta_1$ satisfy $\beta < \beta_1 < 1$,
%% $\gamma < 1$ because $\BE_{\gamma-1}\BT(x)$ is not defined if $\gamma = 1$, because denominator in the obj becomes 0.
then $\BE_{\beta -1}[\BT(x)] \leq \BE_{\beta_1 - 1}[\BT(x)]$.
%
%% \includegraphics[scale = 0.1]{E_a_1_LEQ_E_g_1.JPG}
\end{corollary}
%
\proof{Proof:} The proof follows by replacing $\alpha, \alpha_1$, and $\gamma$ in Proposition \ref{lemma_E_monotonic} with $\beta$, $\beta_1$, and $1$ respectively.
\end{comment}
%
%
%
%
%
%
%
%
%
%\textcolor{red}{Generalization of Lemma 1, given probability levels $\alpha$, $\alpha'$, and $\gamma$ such that $\alpha < \alpha' < \gamma$ then $\BE_{\alpha-\gamma}[\BT(x)] \leq \BE_{\alpha'-\gamma}[\BT(x)]$ }
%
\begin{corollary}\label{prop:monotonicity_nu_star}
For a fixed probability level $\alpha$, if $\gamma_1$ and $\gamma_2$ are such that $\alpha < \gamma_1 \leq \gamma_2$, then $\min_{(x,T(x))\in \XTP} \BE_{\alpha - \gamma_1}\BT(x) \leq \min_{(x,T(x))\in \XTP} \BE_{\alpha - \gamma_2}\BT(x)$. \qed
%For a fixed probability level $\alpha$, $\nu^* = \min_{(x,T(x))\in \XTP} \BE_{\alpha - \gamma}\BT(x)$ as computed in (\ref{eq:bilenear_formulation}) is monotonically increasing in $\gamma\in (\alpha,1]$.
\end{corollary}
The above result follows from Lemma \ref{lemma_E_monotonic}, as for any $(x,T(x))\in \XTP$ and for the choice of $\alpha, \gamma_1, \gamma_2$ satisfying $\alpha < \gamma_1 \leq \gamma_2$, we have $\BE_{\alpha - \gamma_1}\BT(x) \leq  \BE_{\alpha - \gamma_2}\BT(x)$. Thus, for a fixed probability level $\alpha$, $\nu^* = \min_{(x,T(x))\in \XTP} \BE_{\alpha - \gamma}\BT(x)$ as computed in (\ref{eq:bilenear_formulation}) is non-decreasing in $\gamma$.
%
\begin{comment}
\proof
%For a given probability level $\alpha$, we first note that the optimal ($s^*,\theta^*)$ obtained from (\ref{eq:bilenear_formulation}) corresponding to $\nu^*$ remains optimal to (\ref{eq:bilenear_formulation}) irrespective of the value of $\gamma\in (\alpha,1]$. 
%
For any given $x\in \XP$ and known realizations $T(x) = (T_1(x),\dots,T_n(x))$ of $\BT(x)$, we leverage (\ref{eq:Exp_comp_alt_method}), the dual formulation for $\BE_{\alpha - 1}\BT(x)$ in (\ref{eq:gamma_w_0}), and the definition of $-\BE_{\gamma - 1}\BT(x)$ in (\ref{LP_for_E_g1}), to equivalently write the objective function in (\ref{eq:bilenear_formulation}) as: 
%
\begin{align*}
& \frac{1}{\gamma - \alpha}\bigl( (1-\alpha) \BE_{\alpha-1}[\BT(x)] - (1-\gamma) \BE_{\gamma - 1}[\BT(x)]\bigr) 
\\
& = \BE_{\alpha -1}[\BT(x)] - \frac{1-\gamma}{\gamma - \alpha}\bigl(\BE_{\gamma -1}[\BT(x)] - \BE_{\alpha -1}[\BT(x)]\bigr)
\\
& = \BE_{\alpha -1}[\BT(x)] - \Bigl(\frac{1-\alpha}{\gamma - \alpha} - 1\Bigr)\bigl(\BE_{\gamma -1}[\BT(x)] - \BE_{\alpha -1}[\BT(x)]\bigr).
\end{align*}
%
With the Assumption \ref{assumption:alpha_gamma}, we note that $\BE_{\gamma -1}[\BT(x)] - \BE_{\alpha -1}[\BT(x)]$ is positive from Lemma \ref{lemma_E_monotonic}, and $\frac{1-\alpha}{\gamma - \alpha} - 1$ is decreasing for $\gamma \in (\alpha,1]$. So, the objective function and hence $\nu^*$ is monotonically increasing in $\gamma$. 
\endproof
\end{comment}
%  

\subsection{Minimizing value-at-risk ($\Var$)}\label{subsubsection_Var_min}

The proposed formulation for computing $\nu^*$ in (\ref{eq:bilenear_formulation}), also allows solving the highly studied $\Var$-minimization problem (see for instance \cite{kunzi2006computational, sarykalin2008value, v2001value, gaivoronski2005value}). As noted in Propositions \ref{proposition_Min_var_relation} and \ref{prop_choose_epsilon}, for a given $\gamma$, if the choice of $\alpha$ in (\ref{eq:bilenear_formulation}) is such that $\alpha \geq \alpha^*$, where $\alpha^*$ is obtained from Algorithm \ref{Algo_chose_alpha} then, (\ref{eq:bilenear_formulation}) computes $\nu^{*}_{V} = \min_{x, T(x) \in \XTP}\Var_{\gamma}\BT(x)$.
An integer programming formulation (IP) as described in (\ref{eq:IP_Var_form}) has been the standard approach for computing $\nu^{*}_{V}$ (see for instance Section 2 in \cite{feng2015practical}), provided the constant $M$ is chosen large enough to ensure that the constraint $T_i(x)\leq \kappa + M$ is not active when $z_i = 0$ for all $i\in [n]$.
\begin{align}\label{eq:IP_Var_form}
    \text{(IP)}: \min_{z,\kappa, x, T(x)}\Bigl\{ \kappa \bigm| & \sum_{i\in [n]} z_i p_i \geq \gamma,\ T_i(x) \leq \kappa + M(1-z_i),
    \\
    & z\in \{0,1\}^{n},\ (x, T(x)) \in \XTP\Bigr\}.\nonumber
\end{align}
In our evaluation section \ie{} Section \ref{section_case_study}, we will report the efficiency of our bilinear optimization framework by comparing our computational results to the results obtained by the above (IP) formulation.

\subsection{Underestimation for minimizing $\alpha - \gamma$ expectation} \label{subsection_underestimate}

Given the NP-hard nature of computing $\nu^*$ in (\ref{eq:bilenear_formulation}), and $\nu^{*}_{V} = \min_{x, T(x) \in \XTP}\Var_{\gamma}\BT(x)$, we will be interested in estimating $\nu^{*}$ and $\nu^{*}_{V}$. In particular, we will construct their underestimates and overestimates in Sections \ref{subsection_underestimate} and \ref{subsection_overestimate} respectively.
As is seen in (\ref{eq:bilenear_formulation}), computing $\nu^*$ involves solving a nonlinear optimization problem. We construct relaxations for (\ref{eq:bilenear_formulation}) to underestimate $\nu^*$. Let (R) be any convex relaxation of (\ref{eq:bilenear_formulation}), in an extended space $(\psi, \zeta)$, where $\psi = (x, T(x),s,\theta,w')$ such that for each feasible $\bar{\psi}$ to (\ref{eq:bilenear_formulation}), there is a feasible $(\bar{\psi}, \bar{\zeta})$ to (R), such that the objective of (R) evaluates to $\nu^*$ or lower. Leveraging the definitions of $\XTP$, $\WprimeP$, and $\STpolytope(\cdot)$, we construct a candidate for (R) in (\ref{eq:Convex_relax_bilenear_formulation}) using the first-level reformulation linearization technique (RLT) \cite{sherali2013reformulation, sherali1992new}.
\begin{align}\label{eq:Convex_relax_bilenear_formulation}
(\text{R}): \min_{\psi, \zeta} \quad & \frac{1}{\gamma - \alpha}\Bigl((1-\alpha)s + \sum_{i\in [n]}\bigl(\theta_i - T_{ii}^{w}(x) p_i \bigr)\Bigr) 
\\
\text{s.t. } & (s,\theta) \in \STpolytope(x,T(x)), \ w'\in \WprimeP, \ (x,T(x)) \in \XTP, \nonumber
\\
& \bigl(\psi, \zeta\bigr) \in \mathcal{R} \nonumber,
\end{align}
where $\STpolytope(x,T(x))$, $\WprimeP$, and $\XTP$ are as defined in (\ref{eq:STpolytope}), (\ref{LP_for_E_g1}), and Assumption \ref{assumption_XT_polytope} respectively. To obtain the constraints defining $\mathcal{R}$ in (\ref{eq:Convex_relax_bilenear_formulation}), we first post-multiply the constraints defining $\STpolytope(x,T(x))$ and $\XTP$ with the constraints defining $\WprimeP$. 
In Table \ref{Table:RLT_product_constraints} we provide the list of these product constraints. For instance, in Table \ref{Table:RLT_product_constraints}, we obtain $\theta_i w'_j + sp_i w'_j - T_i(x)p_i w'_j \geq 0$ and $\theta_i + sp_i - T_i(x)p_i - (\theta_i w'_j + sp_i w'_j - T_i(x)p_i w'_j) \geq 0$ by multiplying $\theta_i + sp_i - T_i(x)p_i \geq 0$ with $w'_j \geq 0$ and $1-w'_j \geq 0$ respectively.
We then introduce new linearization variables for the nonlinear product terms as described in Table \ref{Table:RLT_linearizing_variable}. The optimization occurs over $(\psi, \zeta)$, where $\zeta$ is the set of all linearization variables defined in Table \ref{Table:RLT_linearizing_variable}.

\begin{table}[htbp]
\centering
\renewcommand{\arraystretch}{1.5}
\setlength{\tabcolsep}{8pt} % Adjust column padding
\begin{tabular}{|c|p{0.65\textwidth}|}
\hline
\multicolumn{2}{|c|}{$\STpolytope(x,T(x)): \theta_i + sp_i - T_i(x)p_i \geq 0\ \forall i \in [n]$} \\ 
\hline
$w'_j \geq 0\ \forall j \in [n]$ & 
$\theta_i w'_j + sp_i w'_j - T_i(x)p_i w'_j \geq 0$ \\
\hline
$1 - w'_j \geq 0\ \forall j \in [n]$ & 
$\theta_i + sp_i - T_i(x)p_i - (\theta_i w'_j + sp_i w'_j - T_i(x)p_i w'_j) \geq 0$ \\
\hline
$\displaystyle\sum_{j\in [n]}w'_jp_j - (1-\gamma) = 0$ & 
$\displaystyle \sum_{j\in [n]} \bigl(\theta_iw'_jp_j + sp_i w'_jp_j - T_i(x)p_i w'_j p_j\bigr)$ \\
& $- (1-\gamma)(\theta_i + sp_i - T_i(x)p_i) = 0$ \\
\hline
\hline
\multicolumn{2}{|c|}{$\STpolytope(x,T(x)): \theta_i \geq 0\ \forall i \in [n]$} \\
\hline
$w'_j \geq 0\ \forall j \in [n]$ & 
$\theta_i w'_j \geq 0$ \\
\hline
$1 - w'_j \geq 0\ \forall j \in [n]$ & 
$\theta_i - \theta_iw'_j \geq 0$ \\
\hline
$\displaystyle\sum_{j\in [n]}w'_jp_j - (1-\gamma) = 0$ & 
$\displaystyle \sum_{j\in [n]} \theta_iw'_jp_j - (1-\gamma)\theta_i = 0$ \\
\hline
\hline
\multicolumn{2}{|c|}{$\XTP: \displaystyle\sum_{i\in [k]}A_{ri}x_i + \sum_{i\in [n]}B_{ri}T_i(x) - c_r \leq 0 \ \forall r \in [l]$} \\
\hline
$w'_j \geq 0\ \forall j \in [n]$ & 
$\displaystyle\sum_{i\in [k]}A_{ri}x_iw'_j + \sum_{i\in [n]}B_{ri}T_i(x)w'_j - c_rw'_j \leq 0$ \\
\hline
$1 - w'_j \geq 0\ \forall j \in [n]$ & 
$\displaystyle\sum_{i\in [k]}A_{ri}x_i + \sum_{i\in [n]}B_{ri}T_i(x) - c_r$ \\
& $- \Bigl(\displaystyle\sum_{i\in [k]}A_{ri}x_iw'_j + \sum_{i\in [n]}B_{ri}T_i(x)w'_j - c_rw'_j\Bigr) \leq 0$ \\
\hline
$\displaystyle\sum_{j\in [n]}w'_jp_j - (1-\gamma) = 0$ & 
$\displaystyle\sum_{i\in [k], j \in [n]}A_{ri}x_iw'_j p_j + \sum_{i,j\in [n]}B_{ri}T_i(x)w'_jp_j - c_r \sum_{j\in [n]}w'_jp_j$ \\
& $ -(1 - \gamma)\Bigl(\displaystyle\sum_{i\in [k]}A_{ri}x_i + \sum_{i\in [n]}B_{ri}T_i(x) - c_r\Bigr) = 0$ \\
\hline
\end{tabular}
\hypersetup{pdfborder={0 0 0}} % Temporarily disable hyperlinks

\caption{First-level RLT: Product constraints to be linearized for defining $\mathcal{R}$}
\hypersetup{pdfborder={0 0 1}} % Re-enable hyperlinks

\label{Table:RLT_product_constraints}
\end{table}

\begin{table}[htbp]
\centering
\renewcommand{\arraystretch}{1.2}
\begin{tabular}{|c|c|c|}
\hline
Linearizing variables & Nonlinear terms & For all indices \\
\hline
$T_{ij}^{w}(x)$ & $T_{i}(x)\times w'_j$ & $i, j\in [n]$ 
\\
$x^{w}_{lj}$ & $x_l \times w'_j$ & $l\in [k]$, $j\in [n]$
\\
$\theta^{w}_{ij}$ & $\theta_i \times w'_j$ & $i, j \in [n]$
\\
$s^{w}_i$ & $s \times w'_i$ & $i \in [n]$
\\
\hline
\end{tabular}
\caption{Linearizing variables for RLT of (\ref{eq:bilenear_formulation})}
\label{Table:RLT_linearizing_variable}
\end{table}
%
%\textcolor{green}{Let the constraints $T_i(x) = f(x,y_i)$ for all $i\in [n]$ and the constraints defining $x\in \XP\subseteq \Re^{k}$ be represented as $\bigl\{\bigl(x,T(x)\bigr) \in \Re^{k+n} : Ax + BT(x) \geq c\bigr\}$ for $A\in \Re^{l \times k}$, $B\in \Re^{l\times n}$, and $c\in \Re^{l}$.}
%
%
%\textcolor{red}{\st{and $\XP_T$ be defined as $\bigl\{\bigl(x,T(x)\bigr) \in \Re^{k+n} : Ax + BT(x) \geq c\bigr\}$ for $A\in \Re^{l \times k}$, $B\in \Re^{l\times n}$, and $c\in \Re^{l}$.}
%
%
%
%
\begin{proposition}\label{prop_RLT_exact} 
For probability levels $\alpha, \gamma$ satisfying Assumption \ref{assumption:alpha_gamma}, if the choice of $\gamma$ is such that $1-\gamma<\min_{i\in [n]}p_i$, then $\WprimeP = \{w'\in [0,1]^{n}: \sum_{i\in [n]} w'_ip_i = 1 - \gamma\}$ in (\ref{LP_for_E_g1}) is a simplex, and $\text{R}^* = \nu^*$, where $\text{R}^*$ represents the optimal value of (R) in (\ref{eq:Convex_relax_bilenear_formulation}) and $\nu^*$ is as in (\ref{eq:bilenear_formulation}). \qed
\end{proposition}
%
\begin{comment}
\proof
We first note that if the choice of $\gamma$ is such that $1-\gamma<\min_{i\in [n]}p_i$, then $\nexists$ $i\in [n]$ such that $w'_i = 1$ and $w'\in \WprimeP$.
%
This is easy to see from the constraints $\sum_{i\in [n]}w'_ip_i = 1-\gamma$ and $w'\geq 0$ in $\WprimeP$, which imply that for all $i\in [n]$, $w'_i\leq \frac{1-\gamma}{p_i}$. Now since the choice of $\gamma$ is such that $1-\gamma<\min_{i\in[n]}p_i$, so $w'_i < 1$ for all $i\in [n]$.
%
Moreover in this case, any solution $\hat{w'}\in \WprimeP$ can be written as $\hat{w'} = \sum_{i\in [n]}\lambda_i\phi_i$ \ie{} as a convex combination of $\{\phi_i\}_{i\in[n]}$, where for all $i\in [n]$, $\lambda_i\geq 0$, $\phi_i = \frac{1-\gamma}{p_i}e_i$ for $e_i$ being the i$^{\text{th}}$ standard basis vector of $\Re^{n}$, and $\sum_{i\in [n]}\lambda_i = 1$.
%
So $\WprimeP$ is a $n-1$ dimensional simplex. Finally, the exactness of the relaxation (R) constructed in (\ref{eq:Convex_relax_bilenear_formulation}) follows from \cite{tawarmalani2024new}.
%
%For $\WprimeP$ to be a $n$-simplex it should be writable as a convex hull of its $n+1$ vertices. Let those vertices be $v_0,\dots,v_n$ then $\WprimeP$ is a $n$-simplex if $\{w'\in [0,1]^{n}:\sum_{i\in [n]}w'_ip_i = 1 - \gamma\} = \{\sum_{i=0}^{n}\lambda_iv_i\mid \sum_{i=0}^{n}\lambda_i = 1, \lambda_i \geq 0\ \forall i = 0,\dots,n\}$. How can $\WprimeP$ become a simplex by increasing $\gamma$. How (R) becomes exact when $\WprimeP$ is a simplex.
\endproof
\end{comment}
%
The above result on the exactness of the first-level RLT relaxation (R) constructed in (\ref{eq:Convex_relax_bilenear_formulation}) follows from \cite{tawarmalani2024new}.
The convex relaxation denoted by (R) in (\ref{eq:Convex_relax_bilenear_formulation}) provides an underestimate of $\nu^*$ in (\ref{eq:bilenear_formulation}). 
Leveraging Proposition \ref{prop_choose_epsilon}, if the choice of $\alpha$ is such that $\gamma > \alpha \geq \alpha^*$ ($\alpha^*$ obtained from Algorithm \ref{Algo_chose_alpha}) in (R), then we get an underestimate to the $\Var$-minimization problem \ie{} $\nu^{*}_{V} = \min_{x, T(x) \in \XTP}\Var_{\gamma}\BT(x)$ as described in Section \ref{subsubsection_Var_min}. However, from our numerical case studies in Section \ref{section_case_study}, we observed that the underestimates for $\nu^*_{V}$ provided by (R) can be further improved by leveraging Lemma \ref{lemma_w_wp}. We explore this improvement in Proposition \ref{prop_improving_RLT}.
\begin{proposition}\label{prop_improving_RLT}
   Consider the first-level RLT relaxation of (\ref{eq:bilenear_formulation}) as described by \text{(R)} in (\ref{eq:Convex_relax_bilenear_formulation}). Given the probability level $\gamma$, if $\alpha$ is chosen satisfying Assumption \ref{assumption:alpha_gamma} and $\alpha \geq \alpha^*$ ($\alpha^*$ obtained from Algorithm \ref{Algo_chose_alpha}), then the underestimate for $\nu^*_{V}$ obtained from (R), can be improved by:
   \begin{enumerate}[label=(\alph*)]
   \item \label{item_1} Replacing $\theta_i \geq \theta^{w}_{ii}$ in (R) with $\theta_i = \theta^{w}_{ii}$ for all $i\in [n]$,  
   \item \label{item_2} Replacing $\theta_i + sp_i - T_i(x) p_i \geq 0$ in \text{(R)} with 
   %\textcolor{red}{$\theta^{w}_{ii} + s^{w}_{i} p_i - T_{ii}^{w}(x) p_i = 0$} 
   $\theta_{i} + s^{w}_{i} p_i - T_{ii}^{w}(x) p_i = 0$ and $s - s^{w}_{i} - T_i(x) + T_{ii}^{w}(x) \geq 0$ for all $i\in[n]$, 
   \end{enumerate}
   where $\theta^{w}$, $s^{w}$, and $T^{w}(x)$ are as described in Table \ref{Table:RLT_linearizing_variable}.
\end{proposition}
\proof
Let ($\RI$) denote the formulation obtained from (R) as defined in (\ref{eq:Convex_relax_bilenear_formulation}) by following the items \ref{item_1} and \ref{item_2} described in Proposition \ref{prop_improving_RLT}.
We will show that feas($\RI$) $\subseteq$ feas(R), where feas($\RI$) and feas(R) represent the feasible regions of ($\RI$) and (R) respectively.
Let $(\bar{\psi}, \bar{\zeta})$ be a feasible solution to ($\RI$), where $\bar{\psi} = (\bar{x}, \bar{T}(\bar{x}), \bar{s},\bar{\theta}, \bar{w}' )$ and $\bar{\zeta} = (\bar{s}^{w}, \bar{\theta}^w, \bar{T}^w(\bar{x}), \bar{x}^w)$. Then, $\bar{\theta}_i + \bar{s}^{w}_{i} p_i - \bar{T}^{w}_{ii}(\bar{x}) p_i = 0$ and $\bar{s} - \bar{s}^{w}_{i} - \bar{T}_i(\bar{x}) + \bar{T}_{ii}^{w}(\bar{x}) \geq 0$ for all $i\in[n]$. Implying:
\begin{align*}
    &\ \bar{s}p_i - \bar{s}^{w}_{i}p_i - \bar{T}_i(\bar{x})p_i + \bar{T}_{ii}^{w}(\bar{x})p_i \geq 0 && \ \forall i \in [n] 
    \\
    \Longleftrightarrow & \ \bar{s}p_i - \bar{s}^{w}_{i}p_i - \bar{T}_i(\bar{x})p_i + \bar{T}_{ii}^{w}(\bar{x})p_i + \bar{\theta}_{i} - \bar{\theta}_{i} \geq 0 && \ \forall i \in [n]
    \\
    \Longleftrightarrow & \ \bar{s}p_i - \bar{T}_i(\bar{x})p_i +  \bar{\theta}_{i}  \geq 0 && \ \forall i \in [n],
\end{align*}
since for all $i\in [n]$, $\ p_i \geq 0$ and $\bar{s}^{w}_{i} p_i + \bar{\theta}_i - \bar{T}^{w}_{ii}(\bar{x}) p_i = 0$.
Also, for all $i\in [n]$, constraints $\bar{\theta}_{i} = \bar{\theta}^{w}_{ii}$ in ($\RI$) trivially satisfy $\bar{\theta}_i \geq \bar{\theta}^{w}_{ii}$ in (R).
So any feasible solution to ($\RI$) is feasible to (R).
Hence, if $\text{R}^*$ and $\RI^*$ represent the optimal values of (R) and ($\RI$) respectively, then  $\text{R}^* \leq \RI^*$.
Now to prove that $\RI^*\leq \nu^*$, we will show that the additional constraints in ($\RI$) are valid constraints.
\noindent
\begin{enumerate}[label=\textbf{B\arabic*},align=left]
    \item \label{eq:item1} $\theta_i = \theta^{w}_{ii}$: We first note that for all $i\in [n]$, $\theta_i w'_i \leq \theta_i w_i = \theta_i$, where the inequality follows due to the multiplication of $w_i \geq w'_i$ (already proved in Lemma \ref{lemma_w_wp}) and $\theta \geq 0$. The equality holds since $\theta_i(1-w_i) = 0$ is true due to the complementary slackness property.
    Thus, $\theta^{w}_{ii} \leq \theta_i$ is true for all $i\in [n]$. We also note that this inequality only requires $\alpha, \gamma$ to satisfy Assumption \ref{assumption:alpha_gamma} in (\ref{eq:bilenear_formulation}).
    For proving $\theta_i \leq \theta^{w}_{ii}$, we make use of the assumption that $\alpha$ satisfies $\alpha \geq \alpha^*$, where $\alpha^*$ is obtained from Algorithm \ref{Algo_chose_alpha}. 
    From the definitions of $\WP$ and $\WprimeP$, we have $w_i, w'_i \in [0,1]$ for all $i\in [n]$.
    %\in \{0,1,e\}$ for all $i\in [n]$, where $e\in (0,1)$.
    %
    Now we look at the three cases:
    \begin{enumerate}[label=(\roman*)]
    \item  $w'_i = 1$: $w'_i = 1\Longrightarrow \theta_i (1-w'_i) = 0$ and hence $\theta_i = \theta^{w}_{ii}$ follows trivially as $(1-w'_i = 0) \times (\theta_i \geq 0) \Longrightarrow \theta_i - \theta^{w}_{ii} = 0$ for all $i\in [n]$.
    \item  $w'_i = e$ for $e\in (0,1)$: If $w'_i \in (0,1)$, then $w_i \in 
           (0,1]$ since $w_i \geq w'_i$ from Lemma \ref{lemma_w_wp}. Moreover, for the choice of $\alpha$, satisfying $\alpha^* \leq \alpha <\gamma$, we observe no jumps in $F^{-1}(p)$ for $p\in [\alpha, \gamma]$ (see Proposition \ref{prop_choose_epsilon}). Thus, $w_i$ can not be one, and we get $w_i \in (0,1)$. This results in $\theta_i = 0$ as $(1-w_i)\theta_i = 0$ is true due to complementary slackness. Hence for all $i\in [n]$, we have $\theta_i(1-w'_i) = 0$ and $\theta_i = \theta^{w}_{ii}$ follows.
    \item $w'_i = 0$: If $w'_i = 0$, then $w_i \in [0,1]$ since $w_i            \geq w'_i$. We note that for the choice of $\alpha$                  satisfying $\alpha \geq 
          \alpha^*$, $w_i$ can not be unity, as there are no jumps in $F^{-1}(p)$ for $p\in [\alpha, \gamma]$. So, if $w_i \in [0,1)$, then due to complementary slackness $(1-w_i)\theta_i = 0 \Longrightarrow \theta_i = 0 \Longrightarrow (1-w'_i)\theta_i = 0$ and hence $\theta_i =      \theta^{w}_{ii}$. 
    \end{enumerate}

    \item $\theta_{i} + s^{w}_{i}p_{i} = T^{w}_{ii}(x)p_i$: 
    We refer to the formulations $\Gamma(x)$ and $\Gamma^{D}(x)$ in (\ref{eq:gamma_w_0}), and observe that for all $i\in [n]$, we have:
    $\theta_i (1-w_i) = 0 \text{ and } w_i (sp_i +\theta_i - T_i(x)p_i) = 0$, which follows due to complementary slackness.
    Moreover, for all $i\in [n]$, $(w_i - w'_i)\bigl(sp_i +\theta_i - T_i(x)p_i\bigr) \geq 0$, which follows from (\ref{eq:STpolytope}) and the fact that $w_i \geq w'_i$ as stated in Lemma \ref{lemma_w_wp}.
    %
    %\\
    %\\
    %We refer to the formulations $\Gamma(x)$ and $\Gamma^{D}(x)$ in (17), and observe that due to complementary slackness, for all $i\in [n]$, we have $\theta_i (1-w_i) = 0 \text{ and } w_i (sp_i +\theta_i - T_i(x)p_i) = 0$.
    %
    %Now \textcolor{red}{$\gamma>\alpha$ together with (\ref{eq:bilenear_formulation}) implies that $w_i \geq w'_i$ for all $i\in[n]$, as $\sum_{i=1}^{n}w_ip_i = 1-\alpha > \sum_{i=1}^{n} w'_ip_i = 1-\gamma$}. 
    %
    %So, for all $i$, $(w_i - w'_i)\bigl(sp_i +\theta_i - T_i(x)p_i\bigr) \geq 0$, which follows from $w_i \geq w'_i$ and (\ref{eq:STpolytope}).
    %
    Now from $w_i (sp_i +\theta_i - T_i(x)p_i) = 0$ and $(w_i - w'_i)\bigl(sp_i +\theta_i - T_i(x)p_i\bigr) \geq 0$, we have $w'_i\bigl(sp_i + \theta_i - T_i(x)p_i\bigr) \leq 0$, where $w'_i$ and $sp_i + \theta_i - T_i(x)p_i$ are non negative from (\ref{LP_for_E_g1}) and (\ref{eq:STpolytope}) respectively, implying that $w'_i\bigl(sp_i + \theta_i - T_i(x)p_i\bigr) = 0$.
    The required equality is obtained by first linearizing $w'_i\bigl(sp_i + \theta_i - T_i(x)p_i\bigr) = 0$ to $s^{w}_{i}p_i + \theta^{w}_{ii} - T^{w}_{ii}(x)p_i = 0$ and then finally leveraging $\theta_i = \theta^{w}_{ii}$ from \ref{eq:item1} to get $s^{w}_{i}p_i + \theta_{i} - T^{w}_{ii}(x)p_i = 0$.
    %
    %\begin{equation}\label{eq:complementary_slack}
    %    \theta_i (1-w_i) = 0 \text{ and } w_i (sp_i +\theta_i - T_i(x)p_i) = 0 \quad \ \forall i \in [n].
    %\end{equation}
    %
    %Now $\gamma>\alpha$ implies that $w_i \geq w'_i$ for all $i$, as $\sum_{i=1}^{n}w_ip_i = 1-\alpha > \sum_{i=1}^{n} w'_ip_i = 1-\gamma$.
    %
    %So, 
    %\begin{equation}\label{eq:product_equation}
    %    (w_i - w'_i)\bigl(sp_i +\theta_i - T_i(x)p_i\bigr) \geq 0 \quad \forall i \in [n],
    %\end{equation}
    %follows from $w_i \geq w'_i$ and (\ref{eq:STpolytope}). 
    %
    %From (\ref{eq:complementary_slack}) and (\ref{eq:product_equation}), we have $w'_i\bigl(sp_i + \theta_i - T_i(x)p_i\bigr) \leq 0$, where $w'_i$ and $sp_i + \theta_i - T_i(x)p_i$ are non negative from (\ref{LP_for_E_g1}) and (\ref{eq:STpolytope}) respectively, implying that $w'_i\bigl(sp_i + \theta_i - T_i(x)p_i\bigr) = 0$.
    % 
    %Now, the required equality is implied by using the linearization variables in $w'_i\bigl(sp_i + \theta_i - T_i(x)p_i\bigr) = 0$ to get $s^{w}_{i}p_i + \theta^{w}_{ii} - T^{w}_{ii}(x)p_i = 0$. \textcolor{red}{No where here I used $\alpha < \gamma$ to be close.}
    %
    \item $s - s^{w}_{i} \geq  T_i(x) - T_{ii}^{w}(x)$: To see the validity of this constraint, we first multiply $\theta_i + sp_i - p_i T_i(x) \geq 0$ in (\ref{eq:STpolytope}) with $1 - w'_i \geq 0$ and then linearize the obtained constraint using Table \ref{Table:RLT_linearizing_variable}, to get $\theta_i - \theta^{w}_{ii} + p_i(s - s^{w}_i) - p_i (T_i(x) - T^{w}_{ii}(x))\geq 0$. Then the required inequality follows since $\theta_i - \theta^{w}_{ii} = 0$ from \ref{eq:item1}. \endproof
\end{enumerate}
\endproof

\subsection{Overestimation for minimizing $\alpha - \gamma$ expectation}\label{subsection_overestimate}

In this section, given probability levels $\alpha, \gamma$ satisfying Assumption \ref{assumption:alpha_gamma}, we will explore ways to overestimate $\nu^* = \min_{x,T(x)\in \XTP}\BE_{\alpha-\gamma}[\BT(x)]$ and $\nu^{*}_{V} = \min_{x, T(x) \in \XTP}\Var_{\gamma}\BT(x)$. Together with the underestimation techniques we discussed in Section \ref{subsection_underestimate}, these overestimating techniques will provide tighter approximations for $\nu^*$ and $\nu^*_{V}$.
Our proposed overestimation approaches for $\nu^*$ and $\nu^*_{V}$ leverage the definition of $\BE_{\alpha - \gamma}[\BT(x)]$ in (\ref{eq:Exp_comp_alt_method}), and the bilinear optimization formulation stated in (\ref{eq:bilenear_formulation}).
Next, in Proposition \ref{prop_delta_overestimate}, we first argue that for any $(x,T(x))\in \XTP$, given probability levels $\alpha, \gamma$ satisfying Assumption \ref{assumption:alpha_gamma}, we have: $\BE_{\alpha-\gamma}[\BT(x)] \leq \BE_{(\alpha + \delta)-(\gamma + \delta')}[\BT(x)]$ for any $\delta, \delta'\geq 0$ such that $\alpha+\delta < \gamma + \delta'$. We exploit this to construct overestimates for $\nu^*$.
\begin{proposition}\label{prop_delta_overestimate}
   Given probability levels $\alpha, \gamma$ that satisfy \ref{assumption:alpha_gamma}. 
   If $\delta$, $\delta'\geq 0$ are chosen so that $\widetilde{\alpha} = \alpha + \delta$, $\widetilde{\gamma} = \gamma + \delta'$ and $\widetilde{\alpha} < \widetilde{\gamma}$ then $\nu^* \leq \min_{x,T(x) \in \XTP}\BE_{\widetilde{\alpha}-\widetilde{\gamma}}[\BT(x)]$.
\end{proposition}
\proof
To prove the result, we show that for any given $(x,T(x))\in \XTP$, $\BE_{\alpha-\gamma}[\BT(x)] \leq \BE_{\widetilde{\alpha}-\widetilde{\gamma}}[\BT(x)]$. 
We consider the following two cases: \textbf{(1)} $\alpha < \gamma \leq \widetilde{\alpha} < \widetilde{\gamma}$: %
Using (\ref{eq:Var_ag_exp1}) from the proof of Proposition \ref{proposition_Min_var_relation}, we get the relation: $\Var_{\alpha}\BT(x) \leq \BE_{\alpha-\gamma}[\BT(x)] \leq \Var_{\gamma}\BT(x) \leq \BE_{\gamma-\widetilde{\alpha}}[\BT(x)] \leq \Var_{\widetilde{\alpha}}\BT(x) \leq \BE_{\widetilde{\alpha}-\widetilde{\gamma}}[\BT(x)] \leq \Var_{\widetilde{\gamma}}\BT(x)$, thus $\BE_{\alpha - \gamma}[\BT(x)] \leq  \BE_{\widetilde{\alpha}-\widetilde{\gamma}}[\BT(x)]$, and hence $\nu^* \leq \min_{x,T(x) \in \XTP}\BE_{\widetilde{\alpha}-\widetilde{\gamma}}[\BT(x)]$.
For the other case, \textbf{(2)} $\alpha \leq \widetilde{\alpha} < \gamma \leq \widetilde{\gamma}$: We again use (\ref{eq:Var_ag_exp1}) to get,
\begin{equation}\label{eq:inequality_relation}
\Var_{\alpha}\BT(x) \leq \BE_{\alpha-\widetilde{\alpha}}[\BT(x)] \leq \Var_{\widetilde{\alpha}}\BT(x) \leq \BE_{\widetilde{\alpha}-\gamma}[\BT(x)] \leq \Var_{\gamma}\BT(x) \leq \BE_{\gamma-\widetilde{\gamma}}[\BT(x)].
\end{equation}
%
%where the first two inequalities follow since $\alpha \leq \widetilde{\alpha}$, the third and forth inequalities because $\widetilde{\alpha} < \gamma$, and the last inequality follows as $\gamma \leq \widetilde{\gamma}$.
%
Now, from (\ref{eq:inequality_relation}), and Lemma \ref{lemma_convex_combin_E}, we can write $\BE_{\alpha-\gamma}[\BT(x)]$ and $\BE_{\widetilde{\alpha}-\widetilde{\gamma}}[\BT(x)]$ 
 as a convex combination of
 $\{\BE_{\alpha-\widetilde{\alpha}}[\BT(x)], \BE_{\widetilde{\alpha}-\gamma}[\BT(x)]\}$ 
 and 
 $\{\BE_{\widetilde{\alpha}- \gamma}[\BT(x)], \BE_{\gamma-\widetilde{\gamma}}[\BT(x)]\}$ respectively.
 Thus, $\BE_{\alpha-\widetilde{\alpha}}[\BT(x)] 
 \leq 
 \BE_{\alpha-\gamma}[\BT(x)]
 \leq
 \BE_{\widetilde{\alpha}-\gamma}[\BT(x)]
 \leq 
 \BE_{\widetilde{\alpha}-\widetilde{\gamma}}[\BT(x)]
 \leq 
 \BE_{\gamma-\widetilde{\gamma}}[\BT(x)]$ and hence $\nu^* \leq \min_{x,T(x)\in \XTP}\BE_{\widetilde{\alpha}-\widetilde{\gamma}}[\BT(x)]$ follows since $\BE_{\alpha - \gamma}[\BT(x)] \leq  \BE_{\widetilde{\alpha}-\widetilde{\gamma}}[\BT(x)]$ for all $(x, T(x))\in \XTP$. 
\endproof
\begin{corollary}\label{corollary_VaR_overestimate}
For a given probability level $\gamma$, if $\nu^{*}_{V} = \min_{x, T(x) \in \XTP}\Var_{\gamma}\BT(x)$, then $\nu^*_{V} \leq \min_{x,T(x) \in \XTP}\BE_{\gamma-\widetilde{\gamma}}[\BT(x)]$, where $\gamma, \widetilde{\gamma}$ satisfy: $\gamma < \widetilde{\gamma}$, $\widetilde{\gamma} = \gamma + \delta$ for $\delta \geq 0$. \qed
%$\tilde{\gamma} = \gamma + \delta$ for $\delta > 0$. \qed
\end{corollary}
\begin{remark}
Proposition \ref{prop_delta_overestimate} can also be interpreted using the integrated quantile function $\IQF(\cdot)$. We observe that for any given $(x,T(x))\in \XTP$, and for the choice of $\alpha, \gamma, \widetilde{\alpha}, \widetilde{\gamma}$, satisfying $\alpha \leq \widetilde{\alpha} < \gamma \leq \widetilde{\gamma}$, we have $\text{S}(\widetilde{\alpha}, \widetilde{\gamma}) \geq \text{S}(\alpha - \gamma)$ (see Figure \ref{fig_E_a_g_aa_gg} for reference).

\begin{figure}[htbp]
\centering
\begin{tikzpicture}[scale=1, domain=0:4]
  % Axes
  \draw[->] (0,0.03) -- (4.8,0.03) node[right] {$p$};
  \draw[->] (0,0.03) -- (0,3.7) node[above] {$\IQF(p)$};

  % IQF curve
  \draw[domain= 0:4.02, samples=100, smooth, variable=\x, black, thick]plot (\x, {0.05*exp(\x)});

  % Dotted lines for key points
  \draw[dotted, line width=0.1mm] (4, 0.03) -- (4, 2.72);
  \draw[dotted, line width=0.1mm] (3.5, 0.03) -- (3.5, 1.6);
  \draw[dotted, line width=0.1mm] (2.9, 0.03) -- (2.9, 0.92);
  \draw[dotted, line width=0.1mm] (1.3, 0.03) -- (1.3, 0.18);
  \draw[dotted, line width=0.1mm] (0, 2.72) -- (4, 2.72);
  \draw[dotted, line width=0.1mm] (0, 1.6) -- (3.5, 1.6);
  \draw[dotted, line width=0.1mm] (0, 0.92) -- (2.9, 0.92);
  \draw[dotted, line width=0.1mm] (0, 0.18) -- (1.3, 0.18);

  % Lines and slopes
  \draw[-, line width=0.4mm, red] (1.3, 0.18) -- (4, 2.72);
  \draw[-, line width=0.3mm, blue] (2.9, 0.92) -- (4, 2.72);
  \draw[-, line width=0.3mm, orange] (1.3, 0.18) -- (3.5, 1.6);

  % Key points and labels
  \node at (1.3, -0.2) {$\alpha$};
  \node at (2.9, -0.2) {$\widetilde{\alpha}$};
  \node at (3.5, -0.2) {$\gamma$};
  \node at (4, -0.2) {$\widetilde{\gamma}$};

  \node at (-0.8, 2.72) {$\IQF(\widetilde{\gamma})$};
  \node at (-0.8, 1.6) {$\IQF(\gamma)$};
  \node at (-0.8, 0.92) {$\IQF(\widetilde{\alpha})$};
  \node at (-0.8, 0.18) {$\IQF(\alpha)$};

  \node at (1.3, 0.18) {\textcolor{red}{$\bullet$}};
  \node at (2.9, 0.92) {\textcolor{red}{$\bullet$}};
  \node at (3.5, 1.6) {\textcolor{red}{$\bullet$}};
  \node at (4, 2.72) {\textcolor{red}{$\bullet$}};

  % Legend
  \node[anchor=west] at (6, 2.8) {\textbf{\rule{1cm}{0.1mm}} $\IQF(\cdot)$};
  \node[anchor=west] at (6, 2.2) {\textcolor{red}{\textbf{\rule{1cm}{0.5mm}}} $\text{S}(\alpha-\widetilde{\gamma}) = \BE_{\alpha - \widetilde{\gamma}} [\BT(x)]$};
  \node[anchor=west] at (6, 1.6) {\textcolor{orange}{\textbf{\rule{1cm}{0.5mm}}} $\text{S}(\alpha-\gamma) = \BE_{\alpha - \gamma}[\BT(x)]$};
  \node[anchor=west] at (6, 1.0) {\textcolor{blue}{\textbf{\rule{1cm}{0.5mm}}} $\text{S}(\widetilde{\alpha} -\widetilde{\gamma}) = \BE_{\widetilde{\alpha} - \widetilde{\gamma}}[\BT(x)]$};

\end{tikzpicture}

\caption{Proposition \ref{prop_delta_overestimate} using $\IQF(\cdot)$: $\text{S}(\alpha - \gamma) \leq \text{S}(\widetilde{\alpha} - \widetilde{\gamma})$}

\label{fig_E_a_g_aa_gg}
\end{figure}
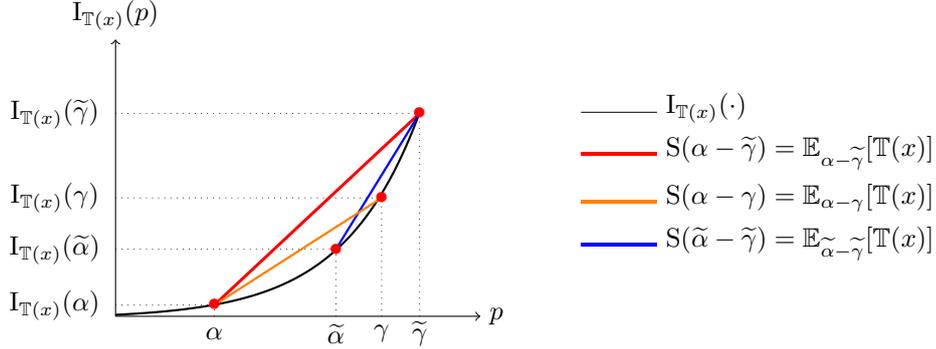
\end{remark}
Leveraging Proposition \ref{prop_delta_overestimate} and Corollary \ref{corollary_VaR_overestimate}, Remark \ref{remark_convex_overestimate} below, describes an approach for obtaining a convex overestimate for $\nu^*$ and $\nu^*_{V}$ by solving the first-level RLT relaxation (R) constructed in (\ref{eq:Convex_relax_bilenear_formulation}).

\begin{remark}\label{remark_convex_overestimate}
Given probability levels $\alpha^*,\gamma^*$ satisfying Assumption \ref{assumption:alpha_gamma}, let $\nu^* = \min_{x,T(x)\in \XTP}\BE_{\alpha^*-\gamma^*}[\BT(x)]$ and $\nu^*_{V} = \min_{x, T(x) \in \XTP}\Var_{\gamma^*}\BT(x)$, then using Propositions \ref{prop_delta_overestimate}, \ref{prop_RLT_exact}, and the RLT relaxation (R) in (\ref{eq:Convex_relax_bilenear_formulation}), a convex overestimate for: 
(1) $\nu^*$ is obtained by solving (R) with $\alpha$ and $\gamma$ being replaced with $\widetilde{\alpha} = \alpha^* + \delta$ and $\widetilde{\gamma} = \gamma^* + \delta'$ respectively where $\delta, \delta' \geq 0$, $\widetilde{\alpha} < \widetilde{\gamma}$, and $\widetilde{\gamma}>1-\min_{i\in[n]}p_i$; 
(2) $\nu^*_{V}$ is obtained by solving (R) with $\alpha$ and $\gamma$ being replaced with $\widetilde{\alpha} = \gamma^*$ and $\widetilde{\gamma} = \gamma^* + \delta'$ respectively where $\widetilde{\alpha}, \widetilde{\gamma}$ satisfy: $\widetilde{\alpha} < \widetilde{\gamma}$, $\widetilde{\gamma} > 1 - \min_{i\in [n]}p_i$ and $\delta' \geq 0$.
%$\delta' \geq 0$ and $\widetilde{\gamma}>1-\min_{i\in[n]}p_i$.
\end{remark}
Remark \ref{remark_convex_overestimate} follows since for this choice of $\widetilde{\gamma}$, $\WprimeP = \{w'\in [0,1]^{n}:\sum_{i\in [n]}w'_ip_i = 1 - \widetilde{\gamma}\}$ is a simplex and hence from Proposition \ref{prop_RLT_exact}, (R) solved for $\widetilde{\alpha}, \widetilde{\gamma}$ is exact to $\min_{x,T(x)\in \XTP}\BE_{\widetilde{\alpha}-\widetilde{\gamma}}[\BT(x)]$. Thus, it provides a valid convex overestimate for $\nu^*$ and $\nu^*_{V}$ from Proposition \ref{prop_delta_overestimate}.
As a special case of Remark \ref{remark_convex_overestimate}, we further note that $\min_{x, T(x)} \BE_{\alpha^*-1}[\BT(x)]$ and $\min_{x, T(x)} \BE_{\gamma^*-1}[\BT(x)]$ are valid convex overestimates for $\nu^* = \min_{x,T(x)\in \XTP}\BE_{\alpha^*-\gamma^*}[\BT(x)]$,  satisfying $\nu^* \leq \min_{x, T(x)} \BE_{\alpha^*-1}[\BT(x)] \leq \min_{x, T(x)} \BE_{\gamma^*-1}[\BT(x)]$. Here, the first and the second inequality hold due to Corollary \ref{prop:monotonicity_nu_star} and Corollary \ref{lemma_E_monotonic} respectively. 

Next, we propose an alternating minimization (AM) approach (see for reference \cite{bezdek2002some, bezdek2003convergence, niesen2007adaptive}) to solve the bilinear optimization formulation in (\ref{eq:bilenear_formulation}). To this, we first define subroutines \Call{Find-Xt}{$\cdot$} and \Call{Find-W}{$\cdot$} such that for a given $w'^* \in \WprimeP$, \Call{Find-Xt}{$w'^*$} solves (\ref{eq:Find-Xt}) and returns its optimal objective value ($\nu_1$) and the optimal solution $(x^*, T^*(x^*))$. Likewise, for a given $(x^*, T^*(x^*)) \in \XTP$, \Call{Find-W}{$x^*, T^*(x^*)$} solves (\ref{eq:Find-W}) and returns its optimal objective value ($\nu_2$) and the optimal solution $w'^*$.
\begin{align}\label{eq:Find-Xt}
\nu_1 = \min_{x,T(x),s,\theta}  \quad &  \frac{1}{\gamma- \alpha} \Bigl((1-\alpha)s + \sum_{i\in [n]}\bigl(\theta_i - T_i(x) p_i w'^{*}_{i}\bigr)\Bigr)
\\
\text{s.t. } &  (s,\theta) \in \STpolytope(x,T(x)), \ (x,T(x))\in \XTP, \nonumber  
\end{align}

\begin{align}\label{eq:Find-W}
\nu_2 = \min_{s,\theta, w'}  \quad &  \frac{1}{\gamma- \alpha} \Bigl((1-\alpha)s + \sum_{i\in [n]}\theta_i - \sum_{i\in [n]}T^{*}_i(x^*) p_i w'_{i}\Bigr)
\\
\text{s.t. } &  (s,\theta) \in \STpolytope(x,T(x)), \ w' \in \WprimeP, \nonumber  
\end{align}
where $\STpolytope(x,T(x))$, $\WprimeP$, and $\XTP$ are as defined in (\ref{eq:STpolytope}), (\ref{LP_for_E_g1}), and Assumption \ref{assumption_XT_polytope}  respectively.
We then leverage subroutines \Call{Find-Xt}{$\cdot$} and \Call{Find-W}{$\cdot$} in Algorithm \ref{Algo:alter_min}, to compute an overestimate for $\min_{(x, T(x)) \in \XTP}\BE_{\alpha-\gamma}[\BT(x)]$.
\begin{algorithm}[htbp]
\caption{Alternate minimization for overestimating $\min_{(x, T(x)) \in \XTP}\BE_{\alpha-\gamma}[\BT(x)]$}\label{Algo:alter_min}
\textbf{Input:} Probability levels $\alpha, \gamma$ satisfying Assumption \ref{assumption:alpha_gamma}.
\\
\textbf{Output:} $\widetilde{\nu}$ satisfying: $\displaystyle \min_{(x,T(x))\in \XTP} \BE_{\alpha - 1}\BT(x) \geq \widetilde{\nu} \geq \min_{(x, T(x)) \in \XTP}\BE_{\alpha-\gamma}[\BT(x)]$
\vspace{1em}
\\
\textbf{Initialize:} $w' = 0, \epsilon = 1e-6$, $i = 0$
\begin{algorithmic}[htbp]
\Procedure{Alt-Min}{}
    \\
    \label{step2-1}\hspace*{0.4cm} $(\nu_i, x^{*},T^{*}(x^*))\leftarrow$ \Call{Find-Xt}{$w'$}
    \\
    \label{step2-2}\hspace*{0.4cm} $(\nu_{i+1}, w'^*) \leftarrow $ \Call{Find-W}{$x^*,T^*(x^*)$}
    \While{$|\nu_{i+1} - \nu_{i}| > \epsilon$} \label{while_loop}
         \State $i \leftarrow i + 2$
         \State $(\nu_i, x^{*},T^{*}(x^*))\leftarrow$ \Call{Find-Xt}{$w'^*$} \label{step2-3}
         \State $(\nu_{i+1}, w'^*) \leftarrow $ \Call{Find-W}{$x^*,T^*(x^*)$} \label{step2-4}
    \EndWhile
    \State $\widetilde{\nu} \leftarrow \nu_{i}$
    \\
    Return $\widetilde{\nu}$.
\EndProcedure
\end{algorithmic}
\end{algorithm}

\begin{proposition}\label{prop_algo_termination}
Algorithm \ref{Algo:alter_min} terminates in finitely many iterations, returning $\widetilde{\nu}$ such that $\min_{(x,T(x))\in \XTP} \BE_{\alpha - 1}\BT(x) \geq \widetilde{\nu} \geq \nu^* = \min_{(x, T(x)) \in \XTP}\BE_{\alpha-\gamma}[\BT(x)]$.
\end{proposition}
\proof
We note that the optimal values $\nu_i$ computed at the Steps \ref{step2-1}, \ref{step2-2}, \ref{step2-3}, and \ref{step2-4} in Algorithm \ref{Algo:alter_min}, are obtained by solving either (\ref{eq:Find-Xt}) or (\ref{eq:Find-W}), which are both constrained versions of the problem $\min_{(x, T(x)) \in \XTP}\BE_{\alpha-\gamma}[\BT(x)]$ in (\ref{eq:bilenear_formulation}). Thus, $\nu_i$ computed throughout the Algorithm \ref{Algo:alter_min} are overestimates of $\nu^*$, and hence $\widetilde{\nu} \geq \nu^*$ also follows. The first inequality in the result \ie{} $\min_{(x,T(x))\in \XTP} \BE_{\alpha - 1}\BT(x) \geq \widetilde{\nu}$ follows from our initialization step in Algorithm \ref{Algo:alter_min} ($w' = 0$) and Proposition \ref{prop:CVar_a_special_case}. 
For a finite $\nu^*$, given the subroutines \Call{Find-Xt}{$\cdot$}, and \Call{Find-W}{$\cdot$}, we note that the sequence of $\{\nu_i\}_{i\geq 0}$ generated in Algorithm \ref{Algo:alter_min} is non-increasing with an upper bound obtained by solving $\min_{x,T(x)\in \XTP}\BE_{\alpha -1}[\BT(x)]$ at its initial step and with a lower bound of $\nu^*$, hence Algorithm \ref{Algo:alter_min} terminates in finitely many iterations of Step \ref{while_loop}.
\endproof
Given probability level $\gamma^*$, we leverage Algorithm \ref{Algo:alter_min}, to construct an overestimate for $\nu^*_V = \min_{(x,T(x))\in \XTP}\Var_{\gamma^*}\BT(x)$, by replacing $\alpha, \gamma$ in Algorithm \ref{Algo:alter_min}, with $\gamma^*$ and $\gamma^* + \delta'$ respectively where $\delta' > 0$. The overestimate $\widetilde{\nu}$ obtained satisfies: $\min_{(x,T(x))\in \XTP}\CVar_{\gamma^*}\BT(x) \geq \widetilde{\nu} \geq \nu^{*}_{V}$. The relation follows from Proposition \ref{prop_algo_termination}.

%The first inequality here follows from Proposition \ref{prop_algo_termination}. The second inequality follows from our initialization step in Algorithm \ref{Algo:alter_min} \ie{}  $w' = 0$ and Proposition \ref{prop:CVar_a_special_case}. We will explore the performance of these overestimation techniques in the evaluation section.

\section{Case study: Minimizing value-at-risk}\label{section_case_study}

In this section, we will explore an application of minimizing value-at-risk. Section \ref{section_NR} will introduce and provide evaluation results for the network traffic engineering problem.

\subsection{Network traffic engineering problem}\label{section_NR}

Consider a graph, $G(V,E)$, where $V$ and $E$ are the set of nodes and edges in $G$ respectively. Let $d:V \times V \to \Re$ represent the traffic demand between node-pairs and $c:E \to\Re$ represent the link capacities. 
For any  source $(s\in V)$ - destination $(t\in V)$ pair, referred to as $\Dp = (s,t)$, we represent the traffic demand between $s$ and $t$ by $d_\Dp$.
The goal of the network traffic engineering problem is to route traffic on the network, obeying all link capacities, with the objective of satisfying all $s-t$ demand requirements.
Frequent network failures in the form of link or node failures pose many challenges to this task \cite{chandra2022probability, chang2019lancet, chandra2024computing, liu2014traffic}.
The availability of empirical data necessary to capture the probability that the network experiences different failures of the link or nodes helps to identify critical failure scenarios \cite{bogle2019teavar}.
Here, a failure scenario refers to the state of the links in the network when a particular set of links or nodes ceases to function.  
For instance, if a failure scenario $(f)$ represents the failure of links in $L_f = \{l_1, l_2\}$, then it represents the network state where all links in the set $E \setminus L_f$ are functional and the links in $L_f$ are non-functional.
Let $Q$ be the set of all such network states, such that $q\in Q$, has an associated probability of $p_q$.
With failures occurring in the network, the entire traffic demand for most $s-t$ pairs can not be fully met. Given a network traffic routing scheme $x$, we let $t_{\Dp q}(x)$ represent the fraction of demand $d_\Dp$ for the $s-t$ pair $\Dp$, that can not be satisfied under the network state $q\in Q$. 
Thus, for a given traffic routing scheme $x$, $(1 - t_{\Dp q}(x))d_\Dp$ is the total demand for the source-destination pair $\Dp$, which can be met under system state $q\in Q$, that occurs with a probability of $p_q$.
For concreteness in Table \ref{table:Traffic_engineering_min_var}, the Columns labeled $q_1,\dots, q_{\lvert Q \rvert}$, refer to the $\lvert Q \rvert$ system states, and the Rows labeled as Pair 1, $\dots$, Pair $n$, represent $n$ source-destination demand pairs. Given a traffic routing scheme $x$, $t_{\Dp q}(x)$ for all $\Dp \times q \in [n]\times [\lvert Q \rvert]$, represents the fraction of the demand out of $d_\Dp$ which is left unmet under the system state $q$ and routing $x$.
Finally, the Row ``S-Loss'' in Table \ref{table:Traffic_engineering_min_var}, computes the maximum loss across all source-destination pairs for a given system state and a given routing scheme. For a system state $q_j$, we represent $\text{S}_{j}\text{-Loss}$ as: $t_j(x) = \max_{i \in [n]} t_{ij}(x)$. 
\begin{table}[htbp]
\centering
\renewcommand{\arraystretch}{1.5}

\begin{tabular}{|c|c|c|c c  c|}
\hline
$\mathbf{d(s,t)}$ & $\mathbf{q_1}$ & $\mathbf{q_2}$ &
$\mathbf{q_3}$ & \dots & $\mathbf{q_{\lvert Q\rvert}}$
\\
\hline

\textbf{Pair 1} &  $t_{11}(\cdot)$ & $t_{12}(\cdot)$ & $t_{13}(\cdot)$  & \dots  & $t_{1 \lvert Q\rvert}(\cdot)$
\\
\hline

\textbf{Pair 2} & $t_{21}(\cdot)$ & $t_{22}(\cdot)$ & $t_{23}(\cdot)$  & \dots  & $t_{2 \lvert Q\rvert}(\cdot)$
\\
\hline

\vdots & \vdots & \vdots & \vdots & \vdots  & \vdots
\\
\hline

\textbf{Pair $\mathbf{n}$} & $t_{n1}(\cdot)$ & $t_{n2}(\cdot)$ & $t_{n3}(\cdot)$  & \dots  & $t_{n \lvert Q\rvert}(\cdot)$
\\
\hline
\hline

\textbf{S-Loss} & $\textbf{S}\mathbf{_1}\textbf{-Loss: }t_{1}(\cdot)$ & $\textbf{S}\mathbf{_2}\textbf{-Loss: }t_{2}(\cdot)$ & $\textbf{S}\mathbf{_3}\textbf{-Loss: }t_{3}(\cdot)$ & \dots  & $\textbf{S}_\mathbf{{\lvert Q \rvert}}\textbf{-Loss: }t_{\lvert Q \rvert}(\cdot)$
\\
\hline
\end{tabular}
 \caption{Loss function for the network traffic engineering problem \cite{rao2021flomore}}
 \label{table:Traffic_engineering_min_var}
\end{table}
\subsubsection{Network traffic engineering - Optimization framework}\label{section_NTE_formulation}

We consider a network topology, represented as a graph $G(V,E)$. 
%\textcolor{red}{For link $e\in E$, $c_e$ is its link capacity.} 
We let $P$ represent the set of all source-destination pairs. Each pair $\Dp\in P$ is associated with a demand $d_\Dp$ that needs to be satisfied. 
$Q$ is the set of system states corresponding to different failure scenarios, such that a system state $q\in Q$ occurs with a probability of $p_q$. 
For each demand pair $\Dp\in P$, the associated demand $d_\Dp$ can only be satisfied by routing traffic on a set of tunnels $R(\Dp)$, which is an input to the problem. 
Let $Y_{tq}$ be defined for all tunnels $t\in R = \cup_{\Dp\in P} R(\Dp)$ and $q\in Q$, such that $Y_{tq} = 1$ if tunnel $t$ is functional in the system state $q$ and $0$ otherwise. 
For all $\Dp\in P$ and $t\in R(\Dp)$, $X_{\Dp t}$ denotes the allocation of traffic demand $d_\Dp$ on tunnel $t \in R(\Dp)$. As shown in Table \ref{table:Traffic_engineering_min_var}, given a traffic routing allocation scheme $x = \{\{X_{\Dp t}\}_{t\in R(\Dp)}\}_{\Dp\in P}$, $t_{\Dp q}(x)$ represents the fraction of lost demand for the source-destination pair $\Dp\in P$, under system state $q\in Q$, and $t_q(x) = \max_{\Dp\in P}t_{\Dp q}(x)$ for all $q\in Q$.

Given a probability level $\gamma$, our goal is to decide the traffic routing allocation scheme $x = \{\{X_{\Dp t}\}_{t\in R(\Dp)}\}_{\Dp\in P}$ such that $\Var_{\gamma}\BT(x)$ is minimized, where $\BT(x)$ represents the random variable having realization: $T_{q}(x) = t_q(x)$ for all $q\in Q$.
To do this, we solve the following:
\begin{equation}\label{eq:Var_evaluate}
\varrho^* = \min_{x, T(x)} \{\Var_{\gamma} \BT(x)\mid (x, T(x))\in \XTP\},
\end{equation}
where $\XTP$ is the set of $\bigl(x = (X_{\Dp t})_{t\in R(\Dp), \Dp\in P}, T(x) = (t_q(x))_{q\in Q}\bigr)$ satisfying (\ref{eq:Network_polytope}). Hereafter, for notational ease, we drop $x$ from $(t_q(x))_{q\in Q}$ when the context is clear.
\begin{subequations}\label{eq:Network_polytope}
\allowdisplaybreaks
    \begin{alignat}{3}
        &&&t_q \geq 1 - \frac{\sum_{t\in R(\Dp)}X_{\Dp t} Y_{tq}}{d_\Dp}  \quad && \forall \Dp \in P, \ \forall q \in Q \label{eq:max_loss_constraint}\\
        &&&\sum_{\Dp\in P} \sum_{t \in R(\Dp):e\in t}X_{\Dp t} \leq c_e && \forall e \in E \label{eq:link_constraint_over_tunnel}\\
        &&& X_{\Dp t} \geq 0 && \forall t \in R(\Dp), \ \forall \Dp \in P \label{eq:non_negative_tunnel}\\
        &&&t_q \in [0,1] && \forall q\in Q \label{eq:t_q_fraction}.
    \end{alignat}
\end{subequations}
For a given $\Dp\in P$, and $q\in Q$, the right hand side of (\ref{eq:max_loss_constraint}) models $t_{\Dp q}$ \ie{} the fraction of traffic demand lost for the demand pair $\Dp$, under scenario $q$. Thus, (\ref{eq:max_loss_constraint}) ensures that $t_q \geq \max_{\Dp\in P}t_{\Dp q}$. 
The left-hand side of constraints in (\ref{eq:link_constraint_over_tunnel}) captures the total traffic routed on all tunnels that constitute link $e$, and we require this to be at most $c_e$ \ie{} the capacity of the link $e\in E$. To model that all tunnels carry non-negative traffic we impose (\ref{eq:non_negative_tunnel}), and since $t_q$ models the fraction of lost traffic demand, hence $t_q\in [0,1]$ in (\ref{eq:t_q_fraction}) is implied. 
Leveraging the (IP) defined in (\ref{eq:IP_Var_form}), and setting $M = 1$ in it, we derive an integer programming formulation to compute $\varrho^*$ in (\ref{eq:IP_network_engg_Var_min}). Here, $\XTP$ is defined as in (\ref{eq:Network_polytope}). 
\begin{align}\label{eq:IP_network_engg_Var_min}
    \varrho^{*} = \min_{z, a, x, T(x)} \Bigl\{ a \bigm| & \ \sum_{q\in Q} z_q p_q \geq \gamma, \ (x, T(x)) \in \XTP,
    \\ 
    & t_q(x) \leq a + 1-z_q,\ z_q \in \{0,1\} \ \forall q \in Q \Bigr\}\nonumber.
\end{align}

Given a specified performance threshold $\widetilde{\varrho}$, network architects frequently seek a traffic routing allocation scheme $x = \{\{X_{\Dp t}\}_{t\in R(\Dp)}\}_{\Dp\in P}$ that ensures the associated value-at-risk metric, $\Var_{\gamma}\BT(x)$, does not exceed $\widetilde{\varrho}$. Achieving this objective by approximating $\varrho^*$ in (\ref{eq:Var_evaluate}), requires computing tighter under- and over-estimates for $\min_{x,T(x)} \{\Var_{\gamma}T(x) \mid (x,T(x))\in \XTP\}$. Moreover, simply ensuring that the overestimate is at most $\widetilde{\varrho}$, provides a quantifiable guarantee of meeting the performance requirement. Next, we apply our estimation techniques to over- and under-estimate $\varrho^{*}$.
%\\
%\\
%Many times, given a desirable performance threshold level $\widetilde{\varrho}$, the network architects desire to obtain a traffic routing allocation scheme $x = \{\{X_{\Dp t}\}_{t\in R(\Dp)}\}_{\Dp\in P}$, such that the associated $\Var_{\gamma}\BT(x) \leq \widetilde{\varrho}$.
%
%Hence, obtaining tight estimates for $\min_{x,T(x)} \{\Var_{\gamma}T(x) \mid x,T(x)\in \XTP\}$ and ensuring that the overestimate is atmost $\widetilde{\varrho}$, provides a guarantee for the performance requirement.

\textbf{Comparing underestimates for $\varrho^*$:}
\begin{enumerate}[label=\textbf{U\arabic*},align=left]
\item \label{U1} One common underestimate for $\varrho^*$ is obtained as the LP relaxation of (\ref{eq:IP_network_engg_Var_min}). The LP relaxation provides a convex underestimate for $\varrho^*$ by solving (\ref{eq:IP_network_engg_Var_min}) after substituting $\bigl\{z_q \in \{0,1\}\bigr\}_{q\in Q}$ with $\bigl\{z_q \in [0,1]\bigr\}_{q\in Q}$.
\item \label{U2} Leveraging Proposition \ref{prop_improving_RLT}, and the first-level RLT relaxation (R) as described in (\ref{eq:Convex_relax_bilenear_formulation}), we construct another convex underestimate for $\varrho^*$, by defining $\XTP$ in (R) to represent the set of all $\bigl(x = (X_{\Dp t})_{t\in R(\Dp), \Dp\in P}, T(x) = (t_q(x))_{q\in Q}\bigr)$ satisfying (\ref{eq:Network_polytope}) and $\STpolytope(x,T(x))$ in (R) to represent $\{(s,\theta) \mid s p_q + \theta_q \geq t_q(x) p_q, \ \theta_q \geq 0 \ \forall q \in Q\}$.
\end{enumerate}

\textbf{Comparing overestimates for $\varrho^*$:} 
\begin{enumerate}[label=\textbf{O\arabic*},align=left]
\item \label{O1} A common approach to overestimate $\varrho^*$ is via minimizing the conditional-value-at-risk \ie{} $\min_{(x,T(x))\in \XTP}\CVar_{\gamma}\BT(x) $(see for instance in \cite{bogle2019teavar, rockafellar2000optimization}). In (\ref{eq:CVAR_Network_engineering}) we provide the formulation to minimize $\CVar_{\gamma}\BT(x)$, which provides an overestimate $\varrho_{\CVar}^{*}$ to $\varrho^*$ \cite{bogle2019teavar, rao2021flomore}.
\begin{align}\label{eq:CVAR_Network_engineering}
 \varrho_{\CVar}^{*} = \min_{\eta, x, T(x), \phi} \bigl\{ \eta + \frac{1}{1-\gamma} \sum_{q\in Q}\phi_q p_q \mid \ &  \phi_q \geq t_q(x) - \eta \ \forall q \in Q, \nonumber
 \\
 & \phi_q \geq 0 \  \forall q \in Q, (x,T(x)) \in \XTP \bigr\}.
\end{align}
\item \label{O2} Another convex overestimate for $\varrho^*$ in (\ref{eq:IP_network_engg_Var_min}) follows from Remark \ref{remark_convex_overestimate}. To obtain this overestimate we solve the first-level RLT relaxation of,
\begin{align}\label{eq:2nd_convex_overestimate}
    \min_{\psi} \frac{1}{\widetilde{\gamma} - \widetilde{\alpha}}\bigl\{(1-\widetilde{\alpha})s + \sum_{q\in Q}(\theta_q - t_{q} w'_q p_q ) \mid  \ &(s,\theta) \in \STpolytope(x,T(x)), \nonumber
    \\
    &w'\in \WprimeP,\ (x,T(x)) \in \XTP\bigr\},
\end{align}
where $\psi = (x,T(x),s,\theta,w')$, $\WprimeP = \{w' \mid \sum_{q\in Q}w'_qp_q = 1 - \widetilde{\gamma}, w'_q \in [0,1] \ \forall q \in Q\}$, $\STpolytope(x,T(x)) = \{(s,\theta) \mid s p_q + \theta_q \geq t_q p_q, \ \theta_q \geq 0 \ \forall q \in Q\}$, and $\XTP$ is the set of all $\bigl(x = (X_{\Dp t})_{t\in R(\Dp), \Dp\in P}, T(x) = (t_q(x))_{q\in Q}\bigr)$ satisfying (\ref{eq:Network_polytope}).
To overestimate $\varrho^{*}$, the choice of $\widetilde{\alpha}, \widetilde{\gamma}$ is made using Remark \ref{remark_convex_overestimate}.
The RLT relaxation is obtained using the steps involved in the construction of (R) in (\ref{eq:Convex_relax_bilenear_formulation}). 
%
%The overestimate $\varrho_{R}^*$ is obtained by solving the constructed RLT relaxation with $\widetilde{\alpha}, \widetilde{\gamma}$ satisfying $\widetilde{\alpha} \geq \alpha^*$, $\widetilde{\gamma} > 1 - \min_{q\in Q}p_q$, and $\widetilde{\gamma} > \gamma$ where $\alpha^*$ is chosen from Algorithm \ref{Algo_chose_alpha}.
%
%We also obtain an overestimate for $\varrho^*$, using Algorithm \ref{Algo_chose_alpha} when $\alpha$ is chosen such that $\alpha \geq \alpha^*$.
%
%
\item \label{O3} The third overestimate for $\varrho^*$ in (\ref{eq:IP_network_engg_Var_min}) is obtained using the Alternating minimization algorithm defined in Algorithm \ref{Algo:alter_min}. To do this, we first leverage Corollary \ref{corollary_VaR_overestimate} to obtain the relation: $\varrho^* \leq \min_{(x,T(x))\in \XTP} \BE_{\gamma - \widetilde{\gamma}}\BT(x)$, where $\widetilde{\gamma} = \gamma + \delta', \gamma < \widetilde{\gamma}$ for $\delta' \geq 0$. 
We then replace the probability levels $\alpha$ and $\gamma$ in Algorithm \ref{Algo:alter_min} with given $\gamma$ and chosen $\widetilde{\gamma}$ respectively, to get $\widetilde{\nu}$ as the final output which satisfies: $\varrho^* \leq \min_{(x,T(x))\in \XTP} \BE_{\gamma - \widetilde{\gamma}}\BT(x) \leq \widetilde{\nu} \leq \min_{(x,T(x))\in \XTP}\CVar_{\gamma}\BT(x)$. In Table \ref{table:min_var_network_engg}, we report the overestimating results obtained by the Algorithm \ref{Algo:alter_min}, for network topologies, by choosing $\delta'$ to be $0.007$ and $0.01$.
\end{enumerate}
Later in this section, we will compare the quality of the under- and over-estimates obtained using the above methods for real network topologies.

\subsubsection{Computational Evaluation on Network traffic engineering}\label{section_network_engg_evaluate}

In this section, we provide numerical results for estimating $\varrho^* = \min_{x, T(x)} \{\Var_{\gamma} \BT(x)\mid x, T(x)\in \XTP\}$, using the approaches discussed in Section \ref{section_NTE_formulation}. Our evaluations will be based on seven network topologies, as listed in Table \ref{table_topology_paper_2}, where $\lvert V\rvert$ and $\lvert E \rvert$ represent the number of vertices and edges in the network.
\begin{table}[htbp]
\centering
\renewcommand{\arraystretch}{1.2}
\setlength{\tabcolsep}{0.65em}
\begin{tabular}{|l l|l|l|l l|l|l|}
\hline
\textbf{Topology} &  & \textbf{$\lvert V \rvert$} & $\lvert E \rvert$ & \textbf{Topology} &  & \textbf{$\lvert V \rvert$} & $\lvert E \rvert$ 
\\ \hline
Quest & \textbf{(Q)}   & 19                & 60        
&
InternetMCI  & \textbf{(M)}     & 18                & 64      
\\
Highwinds & \textbf{(H)}     & 16                & 58      
&
IBM & \textbf{(I)}     & 17                & 46  
\\
B4 & \textbf{(B)}  & 12                & 38   
&
Sprint & \textbf{(S)}   & 10                & 34   
\\
XeeX & \textbf{(X)}  & 22                & 64   
&
 &    &                &   
\\ \hline
\end{tabular}
\caption{Topologies used in evaluation}
\label{table_topology_paper_2}
\end{table}
These topologies have been taken from \cite{knight2011internet} and \cite{kumar2018semi}.
As in \cite{rao2021flomore}, we recursively remove one-degree vertices from the network topologies so that the networks are not disconnected with one link failure scenarios. For each source-destination pair in the network, we construct two physical tunnels that can be used to route the demand. The choice of these tunnels is made such that they are as disjointed as possible, and the preference is given to construct shorter tunnels when multiple choices are available. We use the gravity model \cite{zhang2005network} to generate demand traffic matrices that result in the maximum link utilization (see \cite{chandra2022probability, chandra2024computing, chang2019lancet} for reference) to be in the interval $[0.5, 0.7]$ across the topologies.

%We use gravity model \cite{zhang2005network} to generate demand traffic matrices with the maximum link utilization (see \cite{chandra2022probability, chandra2024computing, chang2019lancet} for definition) ranging in the interval $[0.5, 0.7]$ across the topologies.
%
%

\textbf{Generating system states (Q):} As in \cite{bogle2019teavar, rao2021flomore}, for all networks, we use Weibull distribution to generate the failure probability of each link. The Weibull parameter was chosen so that the median failure probability is approximately 0.001, which is consistent with the empirical data characterizing failures in wide-area-networks \cite{gill2011understanding, markopoulou2008characterization, suchara2011network}. Given a set of link failure probabilities, failure scenarios are sampled based on the probability of their occurrence. We further assume that link failures are independent \cite{gill2011understanding, markopoulou2008characterization}.
We only consider failure scenarios and their equivalent system state that occur with a probability $\geq 0.001$ \cite{gill2011understanding, markopoulou2008characterization, suchara2011network}.

\textbf{Computational setting:} We estimate $\varrho^* = \min_{x, T(x)} \{\Var_{\gamma} \BT(x)\mid (x, T(x)) \in \XTP\}$ for different values of $\gamma$. Since the network topologies are designed to be resilient towards failures (see for instance \cite{menth2005network, wang2010r3}), the network architects desire that the value-at-risk of the lost demand be less than some given threshold value with a high probability. This inspires us to choose $\gamma \geq 0.80$ in Table \ref{table:min_var_network_engg}.
All of our algorithms and formulations were implemented in Python. The LP and IP models were solved using Gurobi 8.0 \cite{optimization2020gurobi}. The CPU used was Intel Xeon E5-2623 @3.00 GHz.
\begin{table}[hbtp]
%\begin{threeparttable}[t]
\centering
%{\setlength{\tabcolsep}{0.01em}
\resizebox{1.0\textwidth}{!}{
\begin{tabular}{|c|c |c c| c c c c|c c c|}
%\hline
%\textbf{T} & $\mathbf{\gamma}$ & \textbf{IP} (T) & \textbf{IP-LP} (L) & \textbf{RLT-}$\mathbf{\Gamma_{1}}$ (L) & $\textbf{\CVar}_{\mathbf{\gamma}}$ (U) & \textbf{AM} (U)
%\\ 
%\hline

\hline
\textbf{T} ($\mathbf{\gamma}$) & \textbf{IP-true} & \ref{U1} & \ref{O1}
 & \ref{U2} & 
\ref{O2} & \textbf{\ref{O3}(0.01)} & \textbf{\ref{O3}(0.007)} & G1 & Our-G & Imp$\%$
\\
\hline

Q (0.99) & 0.47235	&	\textit{0.33381} &	0.59010	&	0.43981	& 0.58535 &	0.59010 & 0.58008 &      0.25629
         &   0.14027
 & 45.27

\\ \hline

%Q (0.98)	& 0.21532	&	\textit{0}	& 0.56231	&	0.13154	& 0.53927			& 0.31135 &				0.24843

%\\ \hline

Q (0.97)	& 0	&	\textit{0}	& 0.47226	&	0	& 0.45444	 & 0.09043 &					0.08920 & 0.47226	& 0.08920	& 81.11

\\ \hline

%Q (0.96) & 0	&	\textit{0}	& 0.36534	&	0 &	0.34906		&				0	& 0

%\\ \hline

Q (0.95) &	0	&	\textit{0} &	0.29403	&	0	& 0.27963		&				0	& 0 &   0.29403	& 0 &	100

\\ \hline

Q (0.90)	& 0	&	\textit{0}	& 0.14702	&	0	& 0.13840		&				0 &	0 &  0.14702	& 0 &	100

\\ \hline

Q (0.85)	&	0	&	\textit{0}	& 0.09801	&	0	& 0.09196	&					0	& 0 & 0.09801	& 0	& 100

\\ \hline

Q (0.80)	&	0	&	\textit{0} &	0.07351	&	0	& 0.06885				&		0 &	0 &  0.07351	& 0	& 100

\\ \hline \hline
M (0.99)	&	0.56512		& \textit{0.49437}	& 0.58946	&	0.56248	& 0.58946		& 0.58946 & 0.58946 &   0.09509	& 0.02697	& 71.63

\\ \hline
%M & 0.98	&	0.46242	&	0.33168	&0.58403	&	0.46050&	0.57737	&	0.50954 &0.47919

%\\ \hline
M (0.97)	&	0.44358	&	\textit{0.18472}	& 0.57773	&	0.37052&	0.56581	&	0.47919 &0.47919 &  0.39301	& 0.10867	& 72.35

\\ \hline
%M & 0.96	&	0.39243	&	0.06857	&0.56176	&	0.23271	&0.55052		&0.47919 &0.47919

%\\ \hline
M (0.95)	&	0.19836	&	\textit{0}	&0.54525	&	0.09566&	0.53597	&0.44172 &	0.44159 &  0.54525	& 0.34593	& 36.55

\\ \hline
M (0.90)	&	0	&	\textit{0} &	0.39324	&	0&	0.38711&	0	& 
 0 &  0.39324	& 0	& 100

\\ \hline
M (0.85)	&	0	&	\textit{0}	&0.26433&		0&	0.25940&	0	& 0 &   0.26433	& 0	& 100

\\ \hline
M (0.80)	&	0 &		\textit{0}	& 0.19825 &		0&	0.19422	&0	&0 &      0.19825	& 0	& 100

\\ \hline \hline

H (0.99)		& 0.32587 &		\textit{0.18975}	& 0.41544	&	0.25125&	0.40396&	0.41544 &	0.37071 &     0.22570	& 0.11946	& 47.07

\\ \hline

%H &0.98	&	0	&	0&	0.40170	&	0	&0.37812	&0.06988 &	0.00907

%\\ \hline

H (0.97)	&	0	&	\textit{0} &	0.31281	&	0&	0.29158	&0&	0 &  0.31281	& 0	& 100

\\ \hline

%H &0.96	&	0	&	0&	0.23805	&	0&	0.21892	&0&	0

%\\ \hline

H (0.95)	&	0	&	\textit{0}	&0.19044	&	0	& 0.17425	& 0	& 0 &   0.19044	& 0	& 100

\\ \hline

H (0.90)		& 0 &		\textit{0}	&0.09522	&	0&	0.08624	& 0 &	0 &   0.09522	& 0	& 100

\\ \hline

H (0.85)	&	0	&	\textit{0} &	0.06348	&	0&	0.05730	&0&	0 &  0.06348	& 0	& 100

\\ \hline

H (0.80)		&0&		\textit{0}	&0.04761&		0	&0.04290&	0&	0 &  0.04761	& 0	& 100

\\ \hline \hline

I (0.99)	&	0.43317	&	\textit{0.29870}	&0.57873	&	0.39322	&0.53706	&	0.57873 &0.57760 &   0.28003	& 0.18438	& 34.16

\\ \hline
%I (0.98)		&0.18818	&	0.05946	&0.52147	&	0.17246	&0.49755		&0.30892 &0.21713

%\\ \hline

I (0.97)		&0	&	\textit{0}	&0.46471	&	0&	0.44780	&0.18290 &	0 &  0.46471	& 0	& 100

\\ \hline

%I (0.96)		&0&		0	&0.42623	&	0&	0.41001	&	0.15156 &0

%\\ \hline

I (0.95)		&0	&	\textit{0}	&0.36153	&	0	&0.34717&	0	& 0 &   0.36153	& 0	& 100

\\ \hline

I (0.90)		&0		& \textit{0}	&0.18182	&	0	& 0.17272	& 0	& 0 &  0.18182	& 0	& 100

\\ \hline

I (0.85)		&0	&	\textit{0}	&0.12121		&0	&0.11472	& 0	& 0 & 0.12121	& 0	& 100

\\ \hline

I (0.80)		&0	&	\textit{0}	&0.09091	&	0	&0.08588	&0	&0 &  0.09091	& 0	& 100

\\ \hline \hline

B (0.99)	&	0.38004		& \textit{0.23903} &	0.46430	&	0.32746&	0.45345&	0.46430 &	0.41483 &    0.22527	& 0.08737	& 61.22

\\ \hline

%B &0.98 & 0.30829	&	0.02446	&0.44354	&	0.30829	&0.41938		&0.30829 &0.30829

%\\ \hline

B (0.97)	&	0.30829	&	\textit{0}	&0.39846	&	0.30829	& 0.38094	&	0.30829 & 0.30829 &  0.39846 & 0	& 100

\\ \hline

%B &0.96	&	0.15726	&	0	&0.37592&		0.02716	&0.36227	&	0.30829 &0.30829

%\\ \hline

B (0.95)	&	0	&	\textit{0} &	0.33645&		0	&0.32427&	0	&0 & 0.33645	& 0	& 100

\\ \hline

B (0.90)	&	0	&	\textit{0}	&0.16824	&	0	&0.16140	&0&	0 &   0.16824	& 0	& 100

\\ \hline

B (0.85)	&	0&		\textit{0}	&0.11216	&	0	&0.10720&	0&	0 &  0.11216	& 0	& 100

\\ \hline

B (0.80)	&	0	&	\textit{0}	&0.08412&		0&	0.08025&	0	&0 &    0.08412	& 0	& 100

\\ \hline \hline

X (0.99)		&0.67883	&	\textit{0.63328} &	0.69608	&	0.67247	&0.69608	&0.69608 &	0.69608 &    0.06279	& 0.02360 &	62.41

\\ \hline

%X (0.98)		&0.57935	&	0.46015	&0.68260	&	0.56260	& 0.66696 	&0.61155 &	0.59765

%\\ \hline

X (0.97)		&0.49032	&	\textit{0.35699}	&0.67427		&0.47431	&0.66310		& 0.58868 & 0.53203 &   0.31728	& 0.05772	& 81.81

\\ \hline

%X &0.96		&0.40300	&	0.28059	&0.66529	&	0.38900	&0.65582		&0.45836 &0.43616

%\\ \hline

X (0.95)		&0.34273	&	\textit{0.20825}	&0.65213		&0.33237	&0.64430	&0.35004 &	0.33746 & 0.44388	& 0.00509	& 98.85

\\ \hline

X (0.90)	&	0.09478	&	\textit{0.06800}	&0.56006	&	0.09478	&0.55516	&	0.33238 &0.33238 &    0.49206	& 0.23760	& 51.71

\\ \hline

X (0.85)		&0.09478	&	\textit{0.01756}	&0.48013		 & 0.09478	& 0.47558		& 0.13168 & 0.13168 &   0.46257	& 0.03690	& 92.02

\\ \hline

X (0.80)		&0.09478	&	\textit{0}	&0.39964	&	0.09478 & 0.39632 &0.13168 &	0.13168	 &  0.39964	& 0.03690	& 90.77

\\ \hline \hline

S (0.99)	&	0.17556	&	\textit{0}	&0.42020	&	0.04450	& 0.33381		& 0.42020 & 0.34891 &  0.42020 &	0.30441	& 27.56

\\ \hline

%S &0.98	&	0	&	0	&0.21463&		0	&0.16941&	0&	0

%\\ \hline

S (0.97)		&0	&	\textit{0} &	0.14335		&0	&0.11055	&0&	0 &    0.14335	& 0	& 100

\\ \hline

%S &0.96	&	0	&	0	&0.10752	&	0	&0.08169&	0	&0

%\\ \hline

S (0.95)	&	0		&\textit{0}&	0.08601	&	0	&0.06478&	0	& 0 &   0.08601	& 0	& 100

\\ \hline

S (0.90)	&	0	&	0	&0.04301	&	0	& 0.03183 &	0	& 0 &     0.04301	& 0	& 100

\\ \hline

S (0.85)	&	0	&	\textit{0} &	0.02867 &		0	& 0.02110	& 0	& 0 &   0.02867	& 0	& 100

\\ \hline

S (0.80)	&	0	&	\textit{0}	& 0.02150 &		0 &	0.01578 &	0	& 0 &  0.02150	& 0	& 100

\\ \hline

\end{tabular}
}

 \caption{Results for estimating $\varrho^*$}

 \label{table:min_var_network_engg}
 \begin{tablenotes}
       \scriptsize 
       \footnotesize
       \item[1]  $\text{Entries in Column ``G1''} = \text{Entries in Column ``O1''} - \text{Entries in Column ``U1''}$. 
       \item[2]  $\text{Entries in ``Our-G''} = \text{Entries in ``O3(0.007)''}- \text{Entries in ``U2''}$. 
       \item[3] Imp$\% = \frac{\text{Entries in ``G1''} - \text{Entries in ``Our-G''}}{\text{Entries in ``G1''}}\times 100$.
    \end{tablenotes}
\end{table}

\textbf{Computational results:} The computational results are detailed in Table \ref{table:min_var_network_engg}. The column labeled ``T ($\gamma$)'', reports the network topology and the probability level $\gamma$ being considered. We use the abbreviations for the network topologies as stated in Table \ref{table_topology_paper_2}. In the Column labeled ``IP-true'', we report the true value of $\varrho^*$, obtained by solving the integer program described in (\ref{eq:IP_network_engg_Var_min}). 
The Column labeled ``U1'' reports the underestimate for $\varrho^*$, obtained as the optimal value of the LP-relaxation of (\ref{eq:IP_network_engg_Var_min}), as described in Item \ref{U1}. Column labeled ``U2'' reports underestimates to $\varrho^*$, obtained by solving the first-level RLT relaxation (R) as described in Item \ref{U2}. We note that for all the test instances, the underestimate obtained by our proposed RLT relaxation in Column ``U2'' is at least as good as the estimates obtained in Column ``U1''.
In the Column ``O1'', we report the overestimate for $\varrho^*$, obtained by solving $\min_{(x,T(x))\in \XTP}\CVar_{\gamma}\BT(x)$ as described in Item \ref{O1}. 
The Column labeled ``O2'' reports the overestimate for $\varrho^*$ obtained using Remark \ref{remark_convex_overestimate} and the first-level RLT relaxation construction described in Item \ref{O2}.
Finally, the overestimates reported in Columns labeled ``O3(0.007)'' and ``O3(0.01)'' are obtained using the Alternating minimization algorithm (see Item \ref{O3} for reference). 
Moreover, the values reported in the parenthesis for Columns ``O3(0.007)" and ``O3(0.01)" represent the choice of $\delta'$ to be used in Item \ref{O3}. We observe that the overestimates obtained for $\delta' = 0.007$ in Column ``O3(0.007)'' are at least as good as the overestimates obtained for $\delta' = 0.01$ in Column ``O3(0.01)''. 
This observation follows from Corollary \ref{prop:monotonicity_nu_star}. 
We note that the estimates for $\varrho^*$ reported in Columns ``U2'', ``O2'', ``O3(0.01)'', and ``O3(0.007)'' are obtained as a consequence of our proposed bilinear optimization program in (\ref{eq:bilenear_formulation}). Moreover, as already shown in Proposition \ref{prop:CVar_a_special_case}, the estimates obtained in Column ``O1'' are obtainable as a special case of our proposed approach. 

To evaluate the computational performance of our approach, we report the estimation gaps obtained as the difference between the over- and under-estimates for $\varrho^*$ (\ie{} estimation gap = Overestimate $-$ Underestimate). Column ``G1'', documents the estimation gap using the well-established over- and under-estimates obtained in Columns ``O1'' and ``U1'' respectively. In Column ``Our-G'', we report the estimation gap obtained from our proposed methods. Since, for all our test instances, the overestimates obtained in Column ``O3(0.007)'' are the best, hence we compute our estimation gap using Columns ``O3(0.007)'' and ``U2''. Finally, in the Column labeled ``Imp$\%$'', we quantify the estimation benefits of our proposed methods by computing the percentage reduction achieved in the estimation gap.

\section{Conclusions}

In this paper, we introduced a new two-stage optimization framework for modeling integrated quantile functions, represented by $\BE_{\alpha - \gamma}[\cdot]$, where $\alpha$ and $\gamma$ are probability levels with $\alpha < \gamma$. The proposed framework formulates the computation of $\nu^* = \min_{x, T(x)} \BE_{\alpha - \gamma} [\BT(x)]$ for a discrete random variable $\BT(x)$ as a bilinear optimization problem, which allows for broader quantile-based risk assessment beyond the traditional Value-at-Risk and Conditional value-at-risk ($\Var$ and $\CVar$) metrics. We examined the properties of this bilinear formulation and established a connection between $\BE_{\alpha - \gamma} [\BT(x)]$ and the integrated quantile function of $\BT(x)$. This approach facilitates the construction of over- and under-estimators for $\BE_{\alpha - \gamma}\BT(x)$, thereby providing tighter estimates for $\nu^*$ and improving the precision of risk quantification. 
Furthermore, we showed that computing $\nu^{*}_{V} = \min_{x,T(x)} \Var_{\gamma} \BT(x)$ is a specific case of solving $\min_{x,T(x)} \BE_{\alpha - \gamma}[\BT(x)]$ as $\alpha$ approaches $\gamma$, which supports the development of tighter estimators for $\nu^{*}_{V}$ in the VaR-minimization context. Additionally, we demonstrated that the common overestimate for $\nu^{*}_{V}$, obtained by solving $\min_{x,T(x)} \CVar_{\gamma} \BT(x)$, can be derived as a special case of our bilinear formulation, which in our case reduces to a linear program.
The practical utility of our framework was also illustrated through its application to network traffic engineering problems. In this context, given a probability level $\gamma$ and a set of probabilistic network failure scenarios, our goal was to minimize the $\gamma^{\text{th}}$ quantile of the loss function, where the loss function captured the maximum unmet demand across all source-destination pairs. Empirical results demonstrated the ability of our approach to produce tighter estimates for the Var minimization problem, thereby also assisting in the validation of the traffic network's performance amidst failures.
For our future work, we aim to expand this framework for broader applications, exploring its potential in other risk-sensitive fields that require resilient decision-making under uncertainty.

\bibliography{sn-bibliography}% common bib file

@article{wozabal2012value,
  title={Value-at-risk optimization using the difference of convex algorithm},
  author={Wozabal, David},
  journal={OR spectrum},
  volume={34},
  number={4},
  pages={861--883},
  year={2012},
  publisher={Springer}
}

@article{pang2004global,
  title={On the global minimization of the value-at-risk},
  author={Pang, Jong-shi and Leyffer, Sven},
  journal={Optimization Methods and Software},
  volume={19},
  number={5},
  pages={611--631},
  year={2004},
  publisher={Taylor \& Francis}
}

@incollection{larsen2002algorithms,
  title={Algorithms for optimization of value-at-risk},
  author={Larsen, Nicklas and Mausser, Helmut and Uryasev, Stanislav},
  booktitle={Financial engineering, E-commerce and supply chain},
  pages={19--46},
  year={2002},
  address = {USA},
  publisher={Springer}
}

@inproceedings{kast1998var,
  title={VaR and optimization},
  author={Kast, R and Luciano, E and Peccati, L},
  booktitle={2nd International Workshop on Preferences and Decisions, Trento, July},
  volume={1},
  number={3},
  pages={1998},
  year={1998}
}

@article{nemirovski2007convex,
  title={Convex approximations of chance constrained programs},
  author={Nemirovski, Arkadi and Shapiro, Alexander},
  journal={SIAM Journal on Optimization},
  volume={17},
  number={4},
  pages={969--996},
  year={2007},
  publisher={SIAM}
}

@article{rockafellar2000optimization,
  title={Optimization of conditional value-at-risk},
  author={Rockafellar, R Tyrrell and Uryasev, Stanislav and others},
  journal={Journal of risk},
  volume={2},
  pages={21--42},
  year={2000}
}

@article{joo1980simple,
  title={A simple proof for von Neumann's minimax theorem},
  author={Jo{\'o}, I},
  journal={Acta Sci. Math},
  volume={42},
  pages={91--94},
  year={1980}
}

@incollection{sarykalin2008value,
  title={Value-at-risk vs. conditional value-at-risk in risk management and optimization},
  author={Sarykalin, Sergey and Serraino, Gaia and Uryasev, Stan},
  booktitle={State-of-the-art decision-making tools in the information-intensive age},
  pages={270--294},
  year={2008},
  address   = {USA},
  publisher={Informs}
}

@article{feng2015practical,
  title={Practical algorithms for value-at-risk portfolio optimization problems},
  author={Feng, Mingbin and W{\"a}chter, Andreas and Staum, Jeremy},
  journal={Quantitative Finance Letters},
  volume={3},
  number={1},
  pages={1--9},
  year={2015},
  publisher={Taylor \& Francis}
}

@article{gushchin2017integrated,
  title={Integrated quantile functions: properties and applications},
  author={Gushchin, Alexander A and Borzykh, Dmitriy A},
  journal={Modern Stochastics: Theory and Applications},
  volume={4},
  number={4},
  pages={285--314},
  year={2017},
  publisher={VTeX: Solutions for Science Publishing}
}

@article{embrechts2013note,
  title={A note on generalized inverses},
  author={Embrechts, Paul and Hofert, Marius},
  journal={Mathematical Methods of Operations Research},
  volume={77},
  number={3},
  pages={423--432},
  year={2013},
  publisher={Springer}
}

@incollection{bogle2019teavar,
  title={TEAVAR: striking the right utilization-availability balance in WAN traffic engineering},
  author={Bogle, Jeremy and Bhatia, Nikhil and Ghobadi, Manya and Menache, Ishai and Bj{\o}rner, Nikolaj and Valadarsky, Asaf and Schapira, Michael},
  booktitle={Proceedings of the ACM Special Interest Group on Data Communication},
  pages={29--43},
  year={2019}
}

@article{rao2021flomore,
  title={FloMore: Meeting bandwidth requirements of flows},
  author={Rao, Sanjay and Tawarmalani, Mohit and others},
  journal={arXiv preprint arXiv:2108.03221},
  year={2021}
}

@article{knight2011internet,
  title={The internet topology zoo},
  author={Knight, Simon and Nguyen, Hung X and Falkner, Nickolas and Bowden, Rhys and Roughan, Matthew},
  journal={IEEE Journal on Selected Areas in Communications},
  volume={29},
  number={9},
  pages={1765--1775},
  year={2011},
  publisher={IEEE}
}

@inproceedings{kumar2018semi,
  title={Semi-Oblivious Traffic Engineering: The Road Not Taken},
  author={Kumar, Praveen and Yuan, Yang and Yu, Chris and Foster, Nate and Kleinberg, Robert and Lapukhov, Petr and Lim, Chiun Lin and Soul{\'e}, Robert},
  booktitle={15th USENIX Symposium on Networked Systems Design and Implementation (NSDI 18)},
  pages={157--170},
  year={2018}
}

@inproceedings{zhang2005network,
  title={Network anomography},
  author={Zhang, Yin and Ge, Zihui and Greenberg, Albert and Roughan, Matthew},
  booktitle={Proceedings of the 5th ACM SIGCOMM conference on Internet Measurement},
  pages={30--30},
  year={2005}
}

@article{optimization2020gurobi,
  title={Gurobi optimizer reference manual, 2020},
  author={Optimization, LLC Gurobi and others},
  journal={URL http://www. gurobi. com},
  volume={12},
  year={2020}
}

@inproceedings{gill2011understanding,
  title={Understanding network failures in data centers: measurement, analysis, and implications},
  author={Gill, Phillipa and Jain, Navendu and Nagappan, Nachiappan},
  booktitle={Proceedings of the ACM SIGCOMM 2011 Conference},
  pages={350--361},
  year={2011}
}

@article{markopoulou2008characterization,
  title={Characterization of failures in an operational IP backbone network},
  author={Markopoulou, Athina and Iannaccone, Gianluca and Bhattacharyya, Supratik and Chuah, Chen-Nee and Ganjali, Yashar and Diot, Christophe},
  journal={IEEE/ACM transactions on networking},
  volume={16},
  number={4},
  pages={749--762},
  year={2008},
  publisher={IEEE}
}

@article{suchara2011network,
  title={Network architecture for joint failure recovery and traffic engineering},
  author={Suchara, Martin and Xu, Dahai and Doverspike, Robert and Johnson, David and Rexford, Jennifer},
  journal={ACM SIGMETRICS Performance Evaluation Review},
  volume={39},
  number={1},
  pages={97--108},
  year={2011},
  publisher={ACM New York, NY, USA}
}

@article{tawarmalani2024new,
  title={New finite relaxation hierarchies for concavo-convex, disjoint bilinear programs, and facial disjunctions},
  author={Tawarmalani, Mohit},
  journal={arXiv preprint arXiv:2405.11068},
  year={2024}
}

@article{Kampke2015,
  title={The Generalized Inverse of Distribution Functions},
  author={K{\"a}mpke, Thomas and Radermacher, Franz Josef and K{\"a}mpke, Thomas and Radermacher, Franz Josef},
  journal={Income Modeling and Balancing: A Rigorous Treatment of Distribution Patterns},
  pages={9--28},
  year={2015},
  publisher={Springer}
}

@article{gaivoronski2005value,
  title={Value-at-risk in portfolio optimization: properties and computational approach},
  author={Gaivoronski, Alexei A and Pflug, Georg},
  journal={Journal of risk},
  volume={7},
  number={2},
  pages={1--31},
  year={2005}
}

@article{v2001value,
  title={Value-at-risk based portfolio optimization},
  author={v Puelz, Amy},
  journal={Stochastic optimization: Algorithms and applications},
  pages={279--302},
  year={2001},
  publisher={Springer}
}

@article{kunzi2006computational,
  title={Computational aspects of minimizing conditional value-at-risk},
  author={K{\"u}nzi-Bay, Alexandra and Mayer, J{\'a}nos},
  journal={Computational Management Science},
  volume={3},
  number={1},
  pages={3--27},
  year={2006},
  publisher={Springer}
}

@Book{sherali2013reformulation,
  author    = {Sherali, Hanif D and Adams, Warren P},
  publisher = {Springer Science \& Business Media},
  title     = {A reformulation-linearization technique for solving discrete and continuous nonconvex problems},
  year      = {2013},
  volume    = {31},
  address   = {USA},
}

@article{sherali1992new,
  title={A new reformulation-linearization technique for bilinear programming problems},
  author={Sherali, Hanif D and Alameddine, Amine},
  journal={Journal of Global optimization},
  volume={2},
  pages={379--410},
  year={1992},
  publisher={Springer}
}

@inproceedings{bezdek2002some,
  title={Some notes on alternating optimization},
  author={Bezdek, James C and Hathaway, Richard J},
  booktitle={Advances in Soft Computing—AFSS 2002: 2002 AFSS International Conference on Fuzzy Systems Calcutta, India, February 3--6, 2002 Proceedings},
  pages={288--300},
  year={2002},
  organization={Springer}
}

@article{bezdek2003convergence,
  title={Convergence of alternating optimization},
  author={Bezdek, James C and Hathaway, Richard J},
  journal={Neural, Parallel \& Scientific Computations},
  volume={11},
  number={4},
  pages={351--368},
  year={2003},
  publisher={Dynamic Publishers, Inc. Atlanta, GA, USA}
}

@inproceedings{niesen2007adaptive,
  title={Adaptive alternating minimization algorithms},
  author={Niesen, Urs and Shah, Devavrat and Wornell, Gregory},
  booktitle={2007 IEEE International Symposium on Information Theory},
  pages={1641--1645},
  year={2007},
  organization={IEEE}
}

@article{chang2019lancet,
  title={Lancet: Better network resilience by designing for pruned failure sets},
  author={Chang, Yiyang and Jiang, Chuan and Chandra, Ashish and Rao, Sanjay and Tawarmalani, Mohit},
  journal={Proceedings of the ACM on Measurement and Analysis of Computing Systems},
  volume={3},
  number={3},
  pages={1--26},
  year={2019},
  publisher={ACM New York, NY, USA}
}

@article{chandra2022probability,
  title={Probability estimation via policy restrictions, convexification, and approximate sampling},
  author={Chandra, Ashish and Tawarmalani, Mohit},
  journal={Mathematical Programming},
  volume={196},
  number={1},
  pages={309--345},
  year={2022},
  publisher={Springer}
}

@inproceedings{chandra2024computing,
  title={Computing Inefficiency of Service Networks under Random Failures},
  author={Chandra, Ashish},
  booktitle={IISE Annual Conference. Proceedings},
  pages={1--6},
  year={2024},
  organization={Institute of Industrial and Systems Engineers (IISE)}
}

@inproceedings{liu2014traffic,
  title={Traffic engineering with forward fault correction},
  author={Liu, Hongqiang Harry and Kandula, Srikanth and Mahajan, Ratul and Zhang, Ming and Gelernter, David},
  booktitle={Proceedings of the 2014 ACM Conference on SIGCOMM},
  pages={527--538},
  year={2014}
}

@inproceedings{wang2010r3,
  title={R3: resilient routing reconfiguration},
  author={Wang, Ye and Wang, Hao and Mahimkar, Ajay and Alimi, Richard and Zhang, Yin and Qiu, Lili and Yang, Yang Richard},
  booktitle={Proceedings of the ACM SIGCOMM 2010 conference},
  pages={291--302},
  year={2010}
}

@inproceedings{menth2005network,
  title={Network resilience through multi-topology routing.},
  author={Menth, Michael and Martin, Ruediger},
  booktitle={DRCN},
  pages={271--277},
  year={2005}
}

@article{dhaene2006risk,
  title={Risk measures and comonotonicity: a review},
  author={Dhaene, Jan and Vanduffel, Steven and Goovaerts, Marc J and Kaas, Rob and Tang, Qihe and Vyncke, David},
  journal={Stochastic models},
  volume={22},
  number={4},
  pages={573--606},
  year={2006},
  publisher={Taylor \& Francis}
}

@article{jiang2024also,
  title={ALSO-X\#: Better convex approximations for distributionally robust chance constrained programs},
  author={Jiang, Nan and Xie, Weijun},
  journal={Mathematical Programming},
  pages={1--64},
  year={2024},
  publisher={Springer}
}

@article{emmer2013best,
  title={What is the best risk measure in practice? A comparison of standard measures},
  author={Emmer, Susanne and Kratz, Marie and Tasche, Dirk},
  journal={arXiv preprint arXiv:1312.1645},
  year={2013}
}

@article{embrechts2018quantile,
  title={Quantile-based risk sharing},
  author={Embrechts, Paul and Liu, Haiyan and Wang, Ruodu},
  journal={Operations Research},
  volume={66},
  number={4},
  pages={936--949},
  year={2018},
  publisher={INFORMS}
}

@article{cont2010robustness,
  title={Robustness and sensitivity analysis of risk measurement procedures},
  author={Cont, Rama and Deguest, Romain and Scandolo, Giacomo},
  journal={Quantitative finance},
  volume={10},
  number={6},
  pages={593--606},
  year={2010},
  publisher={Taylor \& Francis}
}

@article{gambrah2014risk,
  title={Risk measures and portfolio optimization},
  author={Gambrah, Priscilla Serwaa Nkyira and Pirvu, Traian Adrian},
  journal={Journal of Risk and Financial Management},
  volume={7},
  number={3},
  pages={113--129},
  year={2014},
  publisher={MDPI}
}

@article{dixit2020project,
  title={Project portfolio selection and scheduling optimization based on risk measure: a conditional value at risk approach},
  author={Dixit, Vijaya and Tiwari, Manoj Kumar},
  journal={Annals of Operations Research},
  volume={285},
  number={1},
  pages={9--33},
  year={2020},
  publisher={Springer}
}

@article{qazi2023supply,
  title={Supply chain risk network value at risk assessment using Bayesian belief networks and Monte Carlo simulation},
  author={Qazi, Abroon and Simsekler, Mecit Can Emre and Formaneck, Steven},
  journal={Annals of Operations Research},
  volume={322},
  number={1},
  pages={241--272},
  year={2023},
  publisher={Springer}
}

@article{babazadeh2018optimisation,
  title={Optimisation of supply chain networks under uncertainty: conditional value at risk approach},
  author={Babazadeh, Reza and Sabbaghnia, Ali},
  journal={International Journal of Management and Decision Making},
  volume={17},
  number={4},
  pages={488--508},
  year={2018},
  publisher={Inderscience Publishers (IEL)}
}

@article{lotfi2020robust,
  title={A robust optimization model for sustainable and resilient closed-loop supply chain network design considering conditional value at risk},
  author={Lotfi, Reza and Mehrjerdi, Yahia Zare and Pishvaee, Mir Saman and Sadeghieh, Ahmad and Weber, Gerhard-Wilhelm},
  journal={Numerical Algebra, Control and Optimization},
  volume={11},
  number={2},
  pages={221--253},
  year={2020},
  publisher={Numerical Algebra, Control and Optimization}
}

@article{atakan2017minimizing,
  title={Minimizing value-at-risk in single-machine scheduling},
  author={Atakan, Semih and B{\"u}lb{\"u}l, Kerem and Noyan, Nilay},
  journal={Annals of Operations Research},
  volume={248},
  pages={25--73},
  year={2017},
  publisher={Springer}
}

@article{sarin2014minimizing,
  title={Minimizing conditional-value-at-risk for stochastic scheduling problems},
  author={Sarin, Subhash C and Sherali, Hanif D and Liao, Lingrui},
  journal={Journal of Scheduling},
  volume={17},
  pages={5--15},
  year={2014},
  publisher={Springer}
}

@article{thormann2024boosted,
  title={The Boosted Difference of Convex Functions Algorithm for Value-at-Risk Constrained Portfolio Optimization},
  author={Thormann, Marah-Lisanne and Vuong, Phan Tu and Zemkoho, Alain B},
  journal={arXiv preprint arXiv:2402.09194},
  year={2024}
}

@article{wozabal2010difference,
  title={A difference of convex formulation of value-at-risk constrained optimization},
  author={Wozabal, David and Hochreiter, Ronald and Pflug, Georg Ch},
  journal={Optimization},
  volume={59},
  number={3},
  pages={377--400},
  year={2010},
  publisher={Taylor \& Francis}
}

@phdthesis{chandra2022techniques,
  title={Techniques for Uncertainty quantification, Risk minimization, with applications to risk-averse decision making},
  author={Chandra, Ashish},
  year={2022},
  school={Purdue University Graduate School}
}

@article{belles2016use,
  title={The use of flexible quantile-based measures in risk assessment},
  author={Belles-Sampera, Jaume and Guill{\'e}n, Montserrat and Santolino, Miguel},
  journal={Communications in Statistics-Theory and Methods},
  volume={45},
  number={6},
  pages={1670--1681},
  year={2016},
  publisher={Taylor \& Francis}
}

@article{mausser1998beyond,
  title={Beyond VaR: From measuring risk to managing risk},
  author={Mausser, Helmut},
  journal={ALGO research quarterly},
  volume={1},
  number={2},
  pages={5--20},
  year={1998}
}

@article{nemirovski2012safe,
  title={On safe tractable approximations of chance constraints},
  author={Nemirovski, Arkadi},
  journal={European Journal of Operational Research},
  volume={219},
  number={3},
  pages={707--718},
  year={2012},
  publisher={Elsevier}
}

@article{taylor2008estimating,
  title={Estimating value at risk and expected shortfall using expectiles},
  author={Taylor, James W},
  journal={Journal of Financial Econometrics},
  volume={6},
  number={2},
  pages={231--252},
  year={2008},
  publisher={Oxford University Press}
}

@article{hendricks1996evaluation,
  title={Evaluation of value-at-risk models using historical data},
  author={Hendricks, Darryll},
  journal={Economic policy review},
  volume={2},
  number={1},
  year={1996}
}

@article{kuester2006value,
  title={Value-at-risk prediction: A comparison of alternative strategies},
  author={Kuester, Keith and Mittnik, Stefan and Paolella, Marc S},
  journal={Journal of Financial Econometrics},
  volume={4},
  number={1},
  pages={53--89},
  year={2006},
  publisher={Oxford University Press}
}
%% if required, the content of .bbl file can be included here once bbl is generated
%\input sn-article.bbl

\end{document}